\documentclass[11pt]{article}
\usepackage{amscd,amsmath,amssymb}
\usepackage[mathscr]{eucal}

\title{\bf FOCK REPRESENTATIONS AND QUANTUM MATRICES}

\author{D. Shklyarov \and S. Sinel'shchikov \and L. Vaksman}
\date{}

\newtheorem{theorem}{Theorem}[section]
\newtheorem{lemma}[theorem]{Lemma}
\newtheorem{proposition}[theorem]{Proposition}
\newtheorem{corollary}[theorem]{Corollary}

\makeatletter \@addtoreset{equation}{section}\makeatother

\begin{document}
\large
\maketitle

{\small \centerline{Institute for Low Temperature Physics \& Engineering}
\centerline{National Academy of Sciences of Ukraine}

\begin{center}
e-mail: vaksman@ilt.kharkov.ua
\end{center}

\medskip


\medskip

 In this paper we study the Fock representation of a certain
$*$-algebra which appears naturally in the framework of quantum group
theory. It is also a generalization of the twisted CCR-algebra introduced by
W. Pusz and S.~Woronowicz. We prove that the Fock representation is a
faithful irreducible representation of the algebra by bounded operators in a
Hilbert space, and, moreover, it is the only (up to unitary equivalence)
representation possessing these properties.

\medskip

 }

\medskip

\section{Introduction}

This work deals with the $*$-algebras $\mathrm{Pol}(\mathrm{Mat}_{m,n})_q$,
$0<q<1$, $1 \le m \le n$, defined by $q$-analogues of the canonical
commutation relations (see section \ref{smr}). Our main result is that the
Fock representation of every such algebra is its only faithful irreducible
$*$-representation by bounded operators.

Let us explain how the $*$-algebra $\mathrm{Pol}(\mathrm{Mat}_{m,n})_q$
appears in quantum group theory. It was demonstrated by Harish-Chandra that
every Hermitian symmetric spaces of non-compact type admits a standard
embedding into a vector space as a bounded symmetric domain.
A $q$-analogue of the Harish-Chandra embedding was constructed in
\cite{SV2}, along with $q$-analogues of the polynomial algebras on vector
spaces. The approach, used in \cite{SV2} in constructing some quantum
polynomial algebras, is similar to the suggestion of V.~Drinfeld \cite{Dr1}
to construct algebras of functions on quantum groups via duality arguments.
In the special case of the symmetric space $SU_{n,m}/S(U_n \times U_m)$ one
has the $*$-algebra $\mathrm{Pol}(\mathcal{M}at_{m,n})_q$ described in
section \ref{da} of \cite{SV2}. Our interest in the algebra
$\mathrm{Pol}(\mathrm{Mat}_{m,n})_q$ is inspired by the isomorphism
\begin{equation}\label{McM}
\mathrm{Pol}(\mathrm{Mat}_{m,n})_q \simeq
\mathrm{Pol}(\mathcal{M}at_{m,n})_q.
\end{equation}
In the special cases $m=1$, $m=n=2$, the $*$-algebras we are interested in
were considered before by W. Pusz, S. Woronowicz \cite{PuW} and by D.
Proskurin, L. Turowska \cite{PT,Tur}. Under these restrictions our result is
already known. In fact, the principal result of this paper is also valid for
all $*$-algebras introduced in \cite{SV2}. This will be shown in a
subsequent work.

Here is the outline of the present paper. In the first part (sections
\ref{smr} -- \ref{bouT}) we sketch the proof of the main result, theorem
\ref{T}. This part presents proofs of only those statements which do not
involve the theory of quantum universal enveloping algebras \cite{Jant}. The
second part (sections \ref{PBW} -- \ref{da}), after recalling the principal
concepts of that theory, presents a construction of the homomorphism of
$*$-algebras
$$
\mathrm{Pol}(\mathrm{Mat}_{m,n})_q \overset{\displaystyle\sim}{\to}\,
\mathrm{Pol}(\mathcal{M}at_{m,n})_q.
$$
The third part (sections \ref{local} -- \ref{ii}) finishes the proof of all
the statements formulated before and, in particular, establishes the
isomorphism (\ref{McM}).

The classical work of W. Arveson \cite{Ar} initiated the research in
non-commutative complex analysis. The study of non-commutative analogues of
function algebras on {\bf bounded} symmetric domains started in \cite{SV2,
V}. Our goal here is to study representations of
$\mathrm{Pol}(\mathrm{Mat}_{m,n})_q$ by {\bf bounded} operators.

The third author thanks H. P. Jakobsen, K. Schm\"udgen, L. Turowska for
helpful discussions.

\bigskip

\section{Statement of the main results}\label{smr}

In what follows $\mathbb{C}$ will be treated as a ground field. We assume
that all the algebras under consideration are unital and $q \in(0,1)$,
unless the contrary is stated explicitly. Consider the well known algebra
$\mathbb{C}[\mathrm{Mat}_{m,n}]_q$ defined by its generators $z_a^\alpha$,
$\alpha=1,\dots,m$; $a=1,\dots,n$, and the commutation relations
\begin{flalign}
& z_a^\alpha z_b^\beta-qz_b^\beta z_a^\alpha=0, & a=b \quad \& \quad
\alpha<\beta,& \quad \text{or}\quad a<b \quad \& \quad \alpha=\beta,
\label{zaa1}
\\ & z_a^\alpha z_b^\beta-z_b^\beta z_a^\alpha=0,& \alpha<\beta \quad
\&\quad a>b,& \label{zaa2}
\\ & z_a^\alpha z_b^\beta-z_b^\beta z_a^\alpha-(q-q^{-1})z_a^\beta
z_b^\alpha=0,& \alpha<\beta \quad \& \quad a<b. & \label{zaa3}
\end{flalign}
This algebra is a quantum analogue of the polynomial algebra
$\mathbb{C}[\mathrm{Mat}_{m,n}]$ on the space of matrices. It will be
convenient for us to introduce an additional assumption $m \le n$.

The algebra admits a natural gradation given by $\deg z_a^\alpha=1$. By
using the diamond lemma \cite{Berg}, a basis of lexicographically ordered
monomials in the vector space $\mathbb{C}[\mathrm{Mat}_{m,n}]_q$ can be
constructed \cite[p.p. 169 -- 171]{DipDo}. This implies that the dimensions
of the corresponding graded components of $\mathbb{C}[\mathrm{Mat}_{m,n}]_q$
and $\mathbb{C}[\mathrm{Mat}_{m,n}]$ are the same.

In a similar way, introduce the algebra
$\mathbb{C}[\overline{\mathrm{Mat}}_{m,n}]_q$, defined by its generators
$(z_a^\alpha)^*$, $\alpha=1,\dots,m$, $a=1,\ldots,n$, and the relations
\begin{flalign}
& (z_b^\beta)^*(z_a^\alpha)^* -q(z_a^\alpha)^*(z_b^\beta)^*=0, \quad a=b
\quad \& \quad \alpha<\beta, \qquad \text{or} & a<b \quad \& \quad
\alpha=\beta, \label{zaa1*}
\\ & (z_b^\beta)^*(z_a^\alpha)^*-(z_a^\alpha)^*(z_b^\beta)^*=0,&
\alpha<\beta \quad\&\quad a>b,& \label{zaa2*}
\\ & (z_b^\beta)^*(z_a^\alpha)^*-(z_a^\alpha)^*(z_b^\beta)^*-
(q-q^{-1})(z_b^\alpha)^*(z_a^\beta)^*=0,& \alpha<\beta \quad \& \quad a<b. &
\label{zaa3*}
\end{flalign}
A gradation in $\mathbb{C}[\overline{\mathrm{Mat}}_{m,n}]_q$ is given by
$\deg(z_a^\alpha)^*=-1$.

Finally, consider the algebra $\mathrm{Pol}(\mathrm{Mat}_{m,n})_q$ whose
generators are $z_a^\alpha$, $(z_a^\alpha)^*$, $\alpha=1,\dots,m$,
$a=1,\dots,n$, and the list of relations is formed by (\ref{zaa1}) --
(\ref{zaa3*}) and
\begin{flalign}
&(z_b^\beta)^*z_a^\alpha=q^2 \cdot \sum_{a',b'=1}^n
\sum_{\alpha',\beta'=1}^mR_{ba}^{b'a'}R_{\beta \alpha}^{\beta'\alpha'}\cdot
z_{a'}^{\alpha'}(z_{b'}^{\beta'})^*+(1-q^2)\delta_{ab}\delta^{\alpha \beta},
& \label{zaa4}
\end{flalign}
with $\delta_{ab}$, $\delta^{\alpha \beta}$ being the Kronecker symbols, and
$$
R_{ij}^{kl}=
\begin{cases}
q^{-1},& i \ne j \quad \& \quad i=k \quad \& \quad j=l
\\ 1,& i=j=k=l
\\ -(q^{-2}-1), & i=j \quad \& \quad k=l \quad \& \quad l>j
\\ 0,& \text{otherwise}.
\end{cases}
$$
The involution in $\mathrm{Pol}(\mathrm{Mat}_{m,n})_q$ is introduced in an
obvious way: $*:z_a^\alpha \mapsto(z_a^\alpha)^*$.

\begin{proposition}\label{1-1}
The linear map
$$
\mathbb{C}[\mathrm{Mat}_{m,n}]_q \otimes
\mathbb{C}[\overline{\mathrm{Mat}}_{m,n}]_q \simeq
\mathrm{Pol}(\mathrm{Mat}_{m,n})_q,\qquad f \otimes g \mapsto f\cdot g
$$
is one-to-one.
\end{proposition}

\smallskip

{\bf Proof.} See Corollary \ref{corol}.

\medskip

\begin{corollary}
The algebra homomorphisms
\begin{gather*}
\mathbb{C}[\mathrm{Mat}_{m,n}]_q \hookrightarrow
\mathrm{Pol}(\mathrm{Mat}_{m,n})_q,\qquad f \mapsto f \otimes 1,
\\ \mathbb{C}[\overline{\mathrm{Mat}}_{m,n}]_q \hookrightarrow
\mathrm{Pol}(\mathrm{Mat}_{m,n})_q,\qquad f \mapsto 1 \otimes f
\end{gather*}
are embeddings.
\end{corollary}

\medskip

We sketch here an explicit construction for a faithful irreducible
$*$-representation of $\mathrm{Pol}(\mathrm{Mat}_{m,n})_q$.
Consider a $\mathrm{Pol}(\mathrm{Mat}_{m,n})_q$-module $\mathcal{H}$
determined by a single generator $v_0$ and the relations
\begin{equation}\label{zv}
(z_a^\alpha)^*v_0=0,\qquad \alpha=1,\ldots,m,\;a=1,\ldots,n.
\end{equation}
It follows easily from proposition \ref{1-1} that

\begin{corollary}\label{sh0}\hfill
\begin{itemize}
\item[i)] $\mathcal{H}=\mathbb{C}[\mathrm{Mat}_{m,n}]_qv_0$;

\item[ii)] there exists a unique sesquilinear form $(.,.)$ on $\mathcal{H}$
with the following properties:
\begin{enumerate}
\item $(v_0,v_0)=1$; \item $(fu,v)=(u,f^*v)$,\ $f \in
\mathrm{Pol}(\mathrm{Mat}_{m,n})_q$,\ $u,v \in \mathcal{H}$.
\end{enumerate}
\end{itemize}
\end{corollary}

{\sc Remark.} To write the form $(.,.)$ explicitly, we introduce a
bigradation on the vector space $\mathrm{Pol}(\mathrm{Mat}_{m,n})_q$:
\begin{equation}\label{bigrad}
\mathrm{Pol}(\mathrm{Mat}_{m,n})_q \simeq \bigoplus_{i,j=0}^\infty
\mathrm{Pol}(\mathrm{Mat}_{m,n})_{q,i,-j}
\end{equation}
with $\mathrm{Pol}(\mathrm{Mat}_{m,n})_{q,i,-j}=
\mathbb{C}[\mathrm{Mat}_{m,n}]_{q,i}\cdot
\mathbb{C}[\overline{\mathrm{Mat}}_{m,n}]_{q,-j}$, and then define a linear
functional $\omega:~\mathrm{Pol}(\mathrm{Mat}_{m,n})_q \to \mathbb{C}$ which
is just the projection onto the (one-dimensional) homogeneous component
$\mathrm{Pol}(\mathrm{Mat}_{m,n})_{q,0,0}\simeq \mathbb{C}$ parallel to
direct sum of all other homogeneous components. Evidently,
$(fv_0,gv_0)=\omega(g^*f)$, $f,g \in \mathbb{C}[\mathrm{Mat}_{m,n}]_q$.

\medskip

\begin{proposition}\label{sfpd}
The above sesquilinear form $(.,.)$ on $\mathcal{H}$ is positive definite.
\end{proposition}

\smallskip

{\bf Proof.} See section \ref{ue}, in particular, proposition \ref{iso}.
\hfill $\square$

\medskip

$\mathcal{H}$ becomes a pre-Hilbert space, so the notion of boundedness for
linear maps in $\mathcal{H}$ makes sense.

\medskip

\begin{proposition}
For every $f \in \mathrm{Pol}(\mathrm{Mat}_{m,n})_q$, the linear map on
$\mathcal{H}$, given by $T(f):v \mapsto fv$, is bounded.
\end{proposition}

\smallskip

{\bf Proof.} See corollary \ref{bou}. \hfill $\square$

\medskip

Thus, $T$ extends up to a $*$-representation $\overline{T}$ of
$\mathrm{Pol}(\mathrm{Mat}_{m,n})_q$ in the Hilbert space
$\overline{\mathcal{H}}$, a completion of $\mathcal{H}$, by bounded
operators.

\medskip

\begin{theorem}\label{T}
$\overline{T}$ is a faithful irreducible $*$-representation of
$\mathrm{Pol}(\mathrm{Mat}_{m,n})_q$ by bounded operators. A representation
with these properties is unique up to unitary equivalence.
\end{theorem}

\smallskip

{\bf Proof.} See proposition \ref{uniq} (uniqueness), proposition \ref{em}
(faithfulness), and section \ref{ue} (irreducibility). \hfill $\square$

\bigskip

\section{An auxiliary algebra \boldmath
$\mathbb{C}[\widetilde{G}]_q$}\label{qphs}

To prove theorem \ref{T}, we need to introduce an auxiliary $*$-algebra
$\mathbb{C}[\widetilde{G}]_q$. Let $N=m+n$. Consider the Hopf algebra
$\mathbb{C}[SL_N]_q$ introduced in the profound works \cite{Dr1, FRT} which
is defined by its generators $\{t_{ij}\}_{i,j=1,\ldots,N}$ and the relations
\begin{flalign}
& t_{\alpha a}t_{\beta b}-qt_{\beta b}t_{\alpha a}=0, & a=b \quad \& \quad
\alpha<\beta,& \quad \text{or}\quad a<b \quad \& \quad \alpha=\beta,
\label{taa1}
\\ & t_{\alpha a}t_{\beta b}-t_{\beta b}t_{\alpha a}=0,& \alpha<\beta \quad
\&\quad a>b,& \label{taa2}
\\ & t_{\alpha a}t_{\beta b}-t_{\beta b}t_{\alpha a}-(q-q^{-1})t_{\beta a}
t_{\alpha b}=0,& \alpha<\beta \quad \& \quad a<b, & \label{taa3}
\\ & \det \nolimits_q \mathbf{t}=1.\label{taa4}
\end{flalign}
Here $\det_q \mathbf{t}$ is a $q$-determinant of the matrix
$\mathbf{t}=(t_{ij})_{i,j=1,\ldots,N}$:
\begin{equation}\label{qdet}
\det \nolimits_q\mathbf{t}=\sum_{s \in
S_N}(-q)^{l(s)}t_{1s(1)}t_{2s(2)}\ldots t_{Ns(N)},
\end{equation}
with $l(s)=\mathrm{card}\{(i,j)|\;i<j \quad \&\quad s(i)>s(j) \}$. The
comultiplication $\Delta$, the counit $\varepsilon$, and the antipode $S$
are defined as follows:
$$
\Delta(t_{ij})=\sum_kt_{ik}\otimes t_{kj},\qquad
\varepsilon(t_{ij})=\delta_{ij},\qquad S(t_{ij})=(-q)^{i-j}\det \nolimits_q
\mathbf{t}_{ji}
$$
with $\mathbf{t}_{ji}$ being the matrix derived from $\mathbf{t}$ by
discarding its $j$-th row and $i$-th column.

Let
$\mathbb{C}[\widetilde{G}]_q\stackrel{\mathrm{def}}{=}(\mathbb{C}[SL_N]_q,*)$
with the involution $*$ been given by
\begin{equation}\label{*}
t_{ij}^*= \mathrm{sign}\left((i-m-1/2)(n-j+1/2)\right)(-q)^{j-i}\det
\nolimits_q \mathbf{t}_{ij}.
\end{equation}
The involution is well defined, as one can see by comparing it to the
involution on the algebra $\mathbb{C}[SU_N]_q$ of regular functions on the
quantum group $SU_N$ (\ref{star5}).\footnote{It worth noting that
$\mathbb{C}[\widetilde{G}]_q$ is not a Hopf $*$-algebra.}

 Recall a standard notation for $q$-minors of $\mathbf{t}$:
$$
t_{IJ}^{\wedge k}\stackrel{\mathrm{def}}{=}\sum_{s \in
S_k}(-q)^{l(s)}t_{i_1j_{s(1)}}\cdot t_{i_2j_{s(2)}}\cdots t_{i_kj_{s(k)}},
$$
with $I=\{(i_1,i_2,\dots,i_k)|\;1 \le i_1<i_2<\dots<i_k \le N \}$,
$J=\{(j_1,j_2,\dots,j_k)|\;1 \le j_1<j_2<\dots<j_k \le N \}$. Introduce the
elements
$$
t=t_{\{1,2,\dots,m \}\{n+1,n+2,\dots,N \}}^{\wedge m},\qquad x=tt^*.
$$
It follows from the definitions that $t$, $t^*$, and $x$ quasi-commute with
all the generators of $\mathbb{C}[SL_N]_q$. More precisely,

\begin{proposition}\label{xcr}
i) $tt^*=t^*t$;

ii) for every polynomial $f$ of a single indeterminate

\begin{gather}
t_{ij} \cdot f(t)=
\begin{cases}
f(qt)t_{ij},& i \le m \quad \&\quad j \le n,
\\ f(q^{-1}t)t_{ij},& i>m \quad \&\quad j>n,
\\ f(t)t_{ij},& \mathrm{otherwise},
\end{cases}\label{ttij}
\\ t_{ij} \cdot f(t^*)=
\begin{cases}
f(qt^*)t_{ij},& i \le m \quad \&\quad j \le n,
\\ f(q^{-1}t^*)t_{ij},& i>m \quad \&\quad j>n,
\\ f(t^*)t_{ij},& \mathrm{otherwise},
\end{cases}\nonumber
\\
t_{ij} \cdot f(x)=
\begin{cases}\label{xtij}
f(q^2x)t_{ij},& i \le m \quad \&\quad j \le n,
\\ f(q^{-2}x)t_{ij},& i>m \quad \&\quad j>n,
\\ f(x)t_{ij},& \mathrm{otherwise}.
\end{cases}
\end{gather}
\end{proposition}

\medskip

Let $\mathbb{C}[\widetilde{G}]_{q,x}$ be the localization of
$\mathbb{C}[\widetilde{G}]_q$ with respect to the multiplicative set
$x,x^2,x^3,\ldots$. An involution in $\mathbb{C}[\widetilde{G}]_{q,x}$ is
imposed in a natural way: $(x^{-1})^*=x^{-1}$. Of course,
$\mathbb{C}[\widetilde{G}]_q \hookrightarrow
\mathbb{C}[\widetilde{G}]_{q,x}$.\footnote{It is well known that
$\mathbb{C}[SL_N]_q$ is a domain \cite[Lemma 9.1.9]{Jos}.} Note that the
element $t$ is invertible in the algebra $\mathbb{C}[\widetilde{G}]_{q,x}$:
$$t^{-1}=t^*x^{-1}.$$

\begin{proposition}\label{emb}
The map
\begin{equation}\label{emb_}
i:z_a^\alpha \mapsto t^{-1}t_{\{1,2,\dots,m \}J_{a \alpha}}^{\wedge m},
\end{equation}
with $J_{a \alpha}=\{n+1,n+2,\dots,N \}\setminus \{N+1-\alpha \}\cup \{a
\}$, admits a unique extension up to an embedding of $*$-algebras
$i:\mathrm{Pol}(\mathrm{Mat}_{m,n})_q \hookrightarrow
\mathbb{C}[\widetilde{G}]_{q,x}$.
\end{proposition}

\smallskip

{\bf Proof.} See proposition \ref{emb3} and the final remarks in section
\ref{ii}. \hfill $\square$

\medskip

 Introduce the element
$y \in \mathrm{Pol}(\mathrm{Mat}_{m,n})_q$ given by
$$
y=1+\sum_{k=1}^m(-1)^k \sum_{\{J'|\;\mathrm{card}(J')=k
\}}\sum_{\{J''|\;\mathrm{card}(J'')=k \}}z_{\quad J''}^{\wedge
kJ'}\left(z_{\quad J''}^{\wedge kJ'}\right)^*.
$$
(The $q$-minors $z_{\quad J''}^{\wedge kJ'}$ for $\mathbf{z}=\{z_a^\alpha
\}$ are defined in the same way as those for $\mathbf{t}$.)
\smallskip

\begin{proposition}\label{yexp}
\begin{equation}\label{y}i(y)=x^{-1}.\end{equation}
\end{proposition}

\smallskip

{\bf Proof.} See section \ref{pry}.

\medskip

 Proposition
\ref{xcr} and formulae (\ref{emb_}), (\ref{y}) imply

\begin{corollary}\label{zy}
For all $\alpha=1,\dots,m$, $a=1,\dots,n$, one has
$$
z_a^\alpha y=q^{-2}yz_a^\alpha,\qquad (z_a^\alpha)^*y=q^2y(z_a^\alpha)^*.
$$
\end{corollary}

\medskip

This allows one to prove the uniqueness statement of our main theorem
\ref{T}.

\begin{proposition}\label{uniq}
A faithful irreducible $*$-representation of
$\mathrm{Pol}(\mathrm{Mat}_{m,n})_q$ by bounded operators in a Hilbert
space, if exists, is unique up to unitary equivalence.
\end{proposition}

\smallskip

{\bf Proof.} Let $\pi$, $\pi'$ be two faithful irreducible
$*$-representations of $\mathrm{Pol}(\mathrm{Mat}_{m,n})_q$ by bounded
linear operators in the Hilbert spaces $H$, $H'$. In particular, $\pi(y)$
and $\pi'(y)$ are non-zero {\sl bounded} self-adjoint operators. The same
standard argument as in \cite{VaSo1} can be used to prove that the non-zero
spectra of the self-adjoint operators $\pi(y)$, $\pi'(y)$ are discrete.
Consider eigenvectors $v$ of $\pi(y)$ and $v'$ of $\pi'(y)$ with $\|v
\|=\|v'\|=1$, associated to a largest modulus eigenvalue of $\pi(y)$ and
$\pi'(y)$, respectively. By a virtue of corollary \ref{zy},
$\pi((z_a^\alpha)^*)v=0$, $\pi'((z_a^\alpha)^*)v'=0$, $\alpha=1,2,\dots,m$,
$a=1,2,\dots,n$. It is easy to show that the kernels of the linear
functionals $(\pi(f)v,v)$, $(\pi'(f)v',v')$ on
$\mathrm{Pol}(\mathrm{Mat}_{m,n})_q$ are just the same subspace $\bigoplus
\limits_{(j,k)\ne(0,0)}\mathbb{C}[\mathrm{Mat}_{m,n}]_{q,j}
\mathbb{C}[\overline{\mathrm{Mat}}_{m,n}]_{q,k}$. Thus
$(\pi(f)v,v)=(\pi'(f)v',v')$. Hence, by irreducibility, the map $v \mapsto
v'$ admits an extension up to a unitary map which intertwines the
representations $\pi$ and $\pi'$. \hfill $\square$

\bigskip

\section{On a $*$-representation of \boldmath $\mathbb{C}[\widetilde{G}]_q$
in a pre-Hilbert space}

Our purpose is to produce a $*$-representation $\widetilde{\mathcal{T}}$ of
$\mathbb{C}[\widetilde{G}]_q$ {\sl in a pre-Hilbert space} such that
$\widetilde{\mathcal{T}}(x)$ is invertible. The representation is a tensor
product of auxiliary representations $\widetilde{\mathcal{T}}_{(k,k+1)}$ of
$\mathbb{C}[SL_N]_q$. The latter are indexed by the standard generators
$(k,k+1)$ of the symmetric group $S_N$.

To describe those auxiliary representations, we need the homomorphisms:
$$
\psi_{(k,k+1)}:\mathbb{C}[SL_N]_q \to \mathbb{C}[SL_2]_q,\qquad
\psi_{(k,k+1)}(t_{ij})=
\begin{cases}
t_{i-k+1,j-k+1},& i,j \in \{k,k+1 \}
\\ \delta_{ij},& \text{otherwise}
\end{cases}
$$
and the following representation $\pi_+$ of $\mathbb{C}[SL_2]_q$ in a vector
space $\mathscr{L}_+$ with a basis $\{e_j \}_{j \in \mathbb{Z}_+}$:
\begin{align}\label{pi+}
\pi_+(t_{12})e_j &=q^{-j}e_j, & \pi_+(t_{21})e_j &=-q^{-(j+1)}e_j,\nonumber
\\
\\ \pi_+(t_{11})e_j &=e_{j+1}, &
\pi_+(t_{22})e_j&= \begin{cases} (1-q^{-2j})e_{j-1}, & j>0,
\\ 0,& j=0.
\end{cases}\nonumber
\end{align}
It is convenient to equip $\mathscr{L}_+$ with structure of a pre-Hilbert
space as follows
$$
(e_i,e_j)=
\begin{cases}
(q^{-2}-1)(q^{-4}-1)\ldots(q^{-2j}-1)\delta_{ij},& j>0,
\\ \delta_{i0},& j=0.
\end{cases}
$$
The representation $\widetilde{\mathcal{T}}_{(k,k+1)}$ of the algebra
$\mathbb{C}[SL_N]_q$ is given by $\widetilde{\mathcal{T}}_{(k,k+1)}{=}\pi_+
\circ \psi_{(k,k+1)}$.

 We now turn to a construction of the representation $\widetilde{\mathcal{T}}$.
Consider the element
$$
u=
\begin{pmatrix}
1 & 2 & \ldots & n & n+1 & n+2 & \ldots & N
\\ m+1 & m+2& \ldots & N & 1 & 2 & \ldots & m
\end{pmatrix},
$$
of the symmetric group $S_N$. This element is a product of cycles $u=s_m
\cdot s_{m-1}\cdot \ldots \cdot s_2 \cdot s_1$, with
\begin{equation}\label{cycle}
s_i=(i,i+1)\cdot(i+1,i+2)\cdot\ldots\cdot(i+n-1,i+n).
\end{equation}
Fix the reduced expression $u=\sigma_1 \sigma_2 \sigma_3 \ldots
\sigma_{mn}$, which is just concatenation of the reduced expressions for the
$s_i$. For example, in the case $m=2$, $n=3$, one has $u=(3,4,5,1,2)$, and
the above reduced expression acquires the form
$u=(2,3)(3,4)(4,5)(1,2)(2,3)(3,4)$.

Now we are in a position to introduce the desired representation:
$$
\widetilde{\mathcal{T}}=\widetilde{\mathcal{T}}_{\sigma_1}\otimes
\widetilde{\mathcal{T}}_{\sigma_2}\otimes \cdots \otimes
\widetilde{\mathcal{T}}_{\sigma_{mn}}.
$$

\begin{proposition}\label{*sa}\hfill \\
i) $\widetilde{\mathcal{T}}$ is a $*$-representation of
$\mathbb{C}[\widetilde{G}]_q$ in the pre-Hilbert space
$\mathscr{L}=\mathscr{L}_+^{\otimes mn}$.
\begin{flalign*}
&ii) & \widetilde{\mathcal{T}}(x){\bf e}_\Bbbk=q^{-2 \sum \limits_jk_j}{\bf
e}_\Bbbk,&&
\end{flalign*}
with ${\bf e}_\Bbbk=e_{k_1}\otimes e_{k_2}\otimes \ldots \otimes
e_{k_{mn}},\qquad \Bbbk=(k_1,k_2,\ldots,k_{mn})\in \mathbb{Z}_+^{mn}$.
\end{proposition}

\medskip

\begin{corollary}
$\widetilde{\mathcal{T}}(x)$ is invertible.
\end{corollary}

\medskip

{\bf Proof} of proposition \ref{*sa}. The first problem is to prove that
$\widetilde{\mathcal{T}}$ is a $*$-representation. The method we apply is
based on well known results of the theory of compact quantum groups
\cite{KS, CP}. That is why we need an involution $\star$ on
$\mathbb{C}[SL_N]_q$, related to the quantum group $SU_N$ \cite{NYM,VaSo2}.
It is given by
\begin{equation}\label{star5}
t_{ij}^\star=(-q)^{j-i}\det \nolimits_q\mathbf{t}_{ij}.
\end{equation}
It is well known that
$\mathbb{C}[SU_N]_q\stackrel{\mathrm{def}}{=}(\mathbb{C}[SL_N]_q,\star)$ is
a Hopf $*$-algebra. Note that its involution $\star$ is related to the
involution $*$ (\ref{*}) as follows:
\begin{equation}\label{star}
t_{ij}^*=\lambda_1(i)\lambda_2(j)t_{ij}^\star,\qquad i,j=1,\ldots,N,
\end{equation}
with
\begin{equation}\label{l*}
\lambda_1(k)=\mathrm{sign}(k-m-1/2),\qquad
\lambda_2(k)=\mathrm{sign}(n-k+1/2).
\end{equation}
Clearly, a representation $\pi$ of $\mathbb{C}[SL_N]_q$ in a pre-Hilbert
space determines a $*$-representation of $\mathbb{C}[\widetilde{G}]_q$ if
and only if $\pi(t_{ij})^*=\lambda_1(i)\lambda_2(j)\pi(t_{ij}^\star)$ for
all $i,j=1,\ldots,N$.

Let $\Lambda'=(\lambda'(1),\lambda'(2),\ldots,\lambda'(N))$,
$\Lambda''=(\lambda''(1),\lambda''(2),\ldots,\lambda''(N))$ be two sequences
whose entries are $\pm 1$. Suppose we are given a representation $\pi$ of
$\mathbb{C}[SL_N]_q$ in a pre-Hilbert space.

{\sc Definition.} $\pi$ is said to be of type $(\Lambda',\Lambda'')$ if
$$\pi(t_{ij})^*=\lambda'(i)\lambda''(j)\pi(t_{ij}^\star).$$

{\sc Remark.} This definition is well illustrated by the special case
$m=n=1$. One has $t_{11}^*=-t_{11}^\star$, $t_{12}^*=t_{12}^\star$,
$t_{21}^*=t_{21}^\star$, $t_{22}^*=-t_{22}^\star$. It is easy to see that
$\pi_+$ is of type $(\Lambda_1,\Lambda_2)$ with $\Lambda_1=(-1,1)$,
$\Lambda_2=(1,-1)$.

\medskip

\begin{lemma}\label{pp}
Suppose that representations $\pi'$ and $\pi''$ are of types
$(\Lambda',\Lambda'')$ and $(\Lambda'',\Lambda''')$ respectively. Then their
tensor product $\pi=\pi'\otimes \pi''$ is of type $(\Lambda',\Lambda''')$.
\end{lemma}

\smallskip

{\bf Proof.} An application of the relation $(\lambda''(k))^2=1$ and the
fact that the comultiplication $\triangle:\mathbb{C}[SU_N]_q \to
\mathbb{C}[SU_N]_q^{\otimes 2}$ is a homomorphism of $*$-algebras yields
\begin{multline*}
(\pi(t_{ij}))^*=\sum_{k=1}^N \pi'(t_{ik})^*\otimes
\pi''(t_{kj})^*=\lambda'(i)\lambda'''(j)\sum_{k=1}^N(\lambda''(k))^2
\pi'(t_{ik}^\star)\otimes \pi''(t_{kj}^\star)=
\\ =\lambda'(i)\lambda'''(j)\pi(t_{ij}^\star)\qquad \square
\end{multline*}

\medskip


\medskip

Turn back to the proof of proposition \ref{*sa}. Consider the sequence
$\Lambda^{(0)}$, $\Lambda^{(1)}$, $\ldots$, $\Lambda^{(mn)}$, given by
$$
\Lambda^{(j)}=\left(\lambda_1(u_{j}(1)),\lambda_1(u_{j}(2)),\ldots,
\lambda_1(u_{j}(N))\right),
$$
with $u_0=e$, $u_1=\sigma_1$, $u_2=\sigma_1 \cdot \sigma_2$, \ldots,
$u_{mn}=u$. Evidently, $\Lambda^{(0)}=\Lambda_1$,
$\Lambda^{(mn)}=\Lambda_2$.
Observe that if for some $j \in \{1,\ldots,N-1 \}$ the pair
$(\Lambda',\Lambda'')$ possesses the properties: $\lambda'(i)=\lambda''(i)$
for $i \notin \{j,j+1 \}$, $\lambda'(j)=-1$, $\lambda''(j)=1$,
$\lambda'(j+1)=1$, $\lambda''(j+1)=-1$ then by the definitions and the
remark before lemma \ref{pp} the representation $\pi_+\circ \psi_{(j,j+1)}$
of $\mathbb{C}[SL_N]_q$ is of type $(\Lambda',\Lambda'')$. In particular,
the representation $\widetilde{\mathcal{T}}_{\sigma_{j+1}}$ is of type
$(\Lambda^{(j)},\Lambda^{(j+1)})$, and the first statement of proposition
\ref{*sa} follows from lemma \ref{pp}.

Turn to the proof of the second statement of the proposition. Recall the
notation
$$
{\bf e}_\Bbbk=e_{k_1}\otimes e_{k_2}\otimes \ldots \otimes e_{k_{mn}},\qquad
\Bbbk=(k_1,k_2,\ldots,k_{mn})\in \mathbb{Z}_+^{mn}
$$
for the standard orthogonal basis of the pre-Hilbert space $\mathscr{L}$.
Let us first demonstrate the idea of the proof in the special case $m=n=2$.
In this case $x=tt^*$ with $t=t_{13}t_{24}-qt_{14}t_{23}$. It follows from
the definitions that
\begin{align*}
\psi_{(2,3)}\otimes \psi_{(3,4)}\otimes \psi_{(1,2)}\otimes
\psi_{(2,3)}(t_{13}) &=1 \otimes 1 \otimes t_{12}\otimes t_{12},
\\ \psi_{(2,3)}\otimes \psi_{(3,4)}\otimes \psi_{(1,2)}\otimes
\psi_{(2,3)}(t_{24}) &= t_{12}\otimes t_{12}\otimes 1 \otimes 1,
\\ \psi_{(2,3)}\otimes \psi_{(3,4)}\otimes \psi_{(1,2)}\otimes
\psi_{(2,3)}(t_{14}) &=0.
\end{align*}
Thus
$$
\psi_{(2,3)}\otimes \psi_{(3,4)}\otimes \psi_{(1,2)}\otimes
\psi_{(2,3)}(t)=t_{12}\otimes t_{12}\otimes t_{12}\otimes t_{12}.
$$
and hence for all $\Bbbk=(k_1,k_2,k_3,k_4)$
\begin{equation}\label{pix2}
\widetilde{\mathcal{T}}(t){\bf e}_\Bbbk=q^{-(k_1+k_2+k_3+k_4)}{\bf e}_\Bbbk.
\end{equation}

Turn to the case of arbitrary $m,n \in \mathbb{N}$. The following statement
generalizes (\ref{pix2}) and implies the second statement of proposition
\ref{*sa}.

\medskip

\begin{lemma}\label{dlk}
For all $\Bbbk \in \mathbb{Z}_+^{mn}$
\begin{equation}\label{te}
\widetilde{\mathcal{T}}(t){\bf e}_\Bbbk=q^{-\sum \limits_jk_j}{\bf e}_\Bbbk.
\end{equation}
\end{lemma}

\smallskip

{\bf Proof.} By definition $\widetilde{\mathcal{T}}=\pi_+^{\otimes mn}\circ
\Psi$ with $$\Psi=\psi_{\sigma_1}\otimes \psi_{\sigma_2}\otimes \ldots
\otimes \psi_{\sigma_{mn}}:\mathbb{C}[SL_N]_q \to
\mathbb{C}[SL_2]_q^{\otimes mn}.$$ To prove (\ref{te}) it suffices to show
that
\begin{equation}\label{tt12}
  \Psi(t)=\underbrace{t_{12}\otimes\ldots\otimes t_{12}}_{mn}.
\end{equation}
Let us prove the latter equality. Denote by $\Psi_{s_i}$ the homomorphism

$$\psi_{(i,i+1)}\otimes
\psi_{(i+1,i+2)}\otimes\ldots\otimes\psi_{(i+n-1,i+n)}:\mathbb{C}[SL_N]_q
\to \mathbb{C}[SL_2]_q^{\otimes n} $$ associated to the cycle $s_i$ (see
(\ref{cycle})). Evidently, $\Psi=\Psi_{s_m}\otimes \Psi_{s_{m-1}}\otimes
\ldots \otimes \Psi_{s_1}$. The following equalities may be deduced easily
from the definition of comultiplication in $\mathbb{C}[SL_N]_q$:

i) if $k<i$ or $l>i+n$ then
$$\Psi_{s_i}(t_{kl})=\delta_{kl}\cdot\underbrace{1 \otimes \ldots
\otimes 1}_{n};$$

ii) if $k=i$ and $l=i+n$ then
$$\Psi_{s_i}(t_{kl})=\underbrace{t_{12} \otimes \ldots \otimes
t_{12}}_{n}.$$ The latter equalities imply
\begin{equation}\label{psitkl}
\Psi(t_{kl})=
\begin{cases}0,& l>k+n
\\ \underbrace{1 \otimes \ldots \otimes 1}_{(m-k)n}\otimes
\underbrace{t_{12}\otimes \ldots \otimes t_{12}}_n \otimes \underbrace{1
\otimes \ldots \otimes 1}_{(k-1)n},& l=k+n
\end{cases}
\end{equation}
where $k\leq m$. Now (\ref{tt12}) follows from (\ref{psitkl}) and the
definition of $t$.
 \hfill $\square$

\bigskip

Since the operator $\widetilde{\mathcal{T}}(x)$ is invertible, the
representation $\widetilde{\mathcal{T}}$ admits a unique extension onto the
$*$-algebra $\mathbb{C}[\widetilde{G}]_{q,x}$, for which we retain the
notation $\widetilde{\mathcal{T}}$.

\bigskip

\section{Unitary equivalence of \boldmath $T$ and $\mathcal{T}$}\label{ue}

Let $\mathcal{T}=\widetilde{\mathcal{T}}\circ i$ be the $*$-representation
of $\mathrm{Pol}(\mathrm{Mat}_{m,n})_q$ deduced from the embedding of
$*$-algebras $i:\mathrm{Pol}(\mathrm{Mat}_{m,n})_q \to
\mathbb{C}[\widetilde{G}]_{q,x}$ described in section \ref{qphs}. We are
about to produce an isomorphism $\mathscr{J}$ of
$\mathrm{Pol}(\mathrm{Mat}_{m,n})_q$-modules $\mathcal{H}$ and
$\mathscr{L}$. Equip the spaces $\mathcal{H}$, $\mathscr{L}$ and the
algebras $\mathbb{C}[\widetilde{G}]_q$, $\mathrm{Pol}(\mathrm{Mat}_{m,n})_q$
with the gradations:
\begin{gather*}
\deg(t_{ij})=
\begin{cases}1,& i \le m \;\&\;j \le n
\\ -1,& i>m \;\&\;j>n
\\ 0,& \text{otherwise},
\end{cases}
\\ \mathscr{L}=\bigoplus_{j=0}^\infty \mathscr{L}_j,\qquad
\mathscr{L}_j=\left \{v \in
\mathscr{L}|\;\widetilde{\mathcal{T}}(x)v=q^{-2j}v \right \},
\\ \deg(z_a^\alpha)=1,\qquad \deg(z_a^\alpha)^*=-1,
\\ \mathcal{H}=\bigoplus_{j=0}^\infty \mathcal{H}_j,\qquad
\mathcal{H}_j=\mathbb{C}[\mathrm{Mat}_{m,n}]_{q,j}v_0.
\end{gather*}
It follows from proposition \ref{xcr} that $\mathscr{L}$ is a graded
$\mathbb{C}[\widetilde{G}]_q$-module. For any homogeneous vector $v \in
\mathscr{L}$ and all $a=1,\ldots,n$, $\alpha=1,\ldots,m$, one has
$$
\deg(\mathcal{T}(z_a^\alpha)v)=\deg(v)+1,\qquad
\deg(\mathcal{T}(z_a^\alpha)^*v)=\deg(v)-1
$$
by a virtue of propositions \ref{xcr} and \ref{emb}. Hence,
\begin{equation}\label{esp}
(T(\psi)v_0,v_0)=(\mathcal{T}(\psi){\bf e}_0,{\bf e}_0),\qquad \psi \in
\mathrm{Pol}(\mathrm{Mat}_{m,n})_q.
\end{equation}
Also, since the vector $v_0$ is cyclic, it follows that the map
$$
\mathscr{J}:\mathcal{H}\to \mathscr{L},\qquad \mathscr{J}:T(\psi)v_0 \mapsto
\mathcal{T}(\psi){\bf e}_0,\qquad \psi \in \mathbb{C}[\mathrm{Mat}_{m,n}]_q,
$$
is a well defined morphism of $\mathrm{Pol}(\mathrm{Mat}_{m,n})_q$-modules
due to proposition \ref{1-1}.

\medskip

\begin{proposition}\label{iso}
The above map $\mathscr{J}:\mathcal{H}\to \mathscr{L}$ is an isomorphism of
$\mathrm{Pol}(\mathrm{Mat}_{m,n})_q$-modules.
\end{proposition}

\medskip

{\sc Remark.} This certainly implies proposition \ref{sfpd} and the unitary
equivalence of representations $T$ and $\mathcal{T}$ (see (\ref{esp})).

\medskip

To prove proposition \ref{iso}, we need several lemmas.

\medskip

\begin{lemma}\label{Izm}
Let $1 \le \alpha_1<\alpha_2<\dots<\alpha_k \le m$, $1 \le a_1<a_2<\dots<a_k
\le n$, $J=\{n+1,n+2,\dots,N \}\setminus
\{n+\alpha_1,n+\alpha_2,\dots,n+\alpha_k \}\cup \{a_1,a_2,\dots,a_k \}$.
Then
$$
i:z_{\quad \{a_1,a_2,\dots,a_k \}}^{\wedge k
\{m+1-\alpha_k,m+1-\alpha_{k-1},\dots,m+1-\alpha_1 \}}\mapsto
t^{-1}t_{\{1,2,\dots,m \}J}^{\wedge m}
$$
\end{lemma}

\smallskip

{\bf Proof.} See section \ref{pry}.

\medskip

\begin{lemma}\label{c=1}
If $v$ is a vector in the space of a $*$-representation $\rho$ of
$\mathbb{C}[\widetilde{G}]_q$ and
\begin{gather}
\rho \left(t_{\{m+1,m+2,\ldots,N \}J}^{\wedge n}\right)v=0,\qquad J \ne
\{1,2,\ldots,n \}\label{pr}
\end{gather}
then $\rho(x)v=v.$
\end{lemma}

\smallskip

{\bf Proof.} The following relation is a consequence of a more general
formula (6.2) of \cite{NYM}. For all $k=1,\ldots,N-1$,
\begin{equation}\label{ws}
\left({t}_{\{1,\ldots,k \}\{N-k+1,\ldots,N \}}^{\wedge k}
\right)^\star=(-q)^{k(N-k)}\cdot {t}_{\{k+1,\ldots,N \}\{1,\ldots,N-k
\}}^{\wedge(N-k)}.
\end{equation}
Recall that $x=tt^*$, with $t=t_{\{1,2,\ldots,m \}\{n+1,n+2,\ldots,N
\}}^{\wedge m}$; it also follows from (\ref{star}), (\ref{ws}) that
$t^*=(-q)^{mn}t_{\{m+1,m+2,\ldots,N \}\{1,2,\ldots,n \}}^{\wedge n}$. Apply
(\ref{pr}) and $\det_q \mathbf{t}=1$ to obtain
$$
v=\rho(\det \nolimits_q \mathbf{t})v=\left((-q)^{mn}t_{\{1,2,\ldots,m
\}\{n+1,n+2,\ldots,N \}}^{\wedge m} \cdot t _{\{m+1,m+2,\ldots,N
\}\{1,2,\ldots,n \}}^{\wedge n}\right)v=\rho(x)v,
$$
which is just our statement. \hfill $\square$

\medskip

The next lemma involves

{\sc Definition.} Let $\mathrm{Pol}(Y)_q$ be a subalgebra of
$\mathbb{C}[\widetilde{G}]_q$ generated by $t_{\{1,2,\dots,m \}I}^{\wedge
m}$, $t_{\{m+1,m+2,\dots,N \}J}^{\wedge n}$, with $\mathrm{card}\,I=m$,
$\mathrm{card}\,J=n$.

It follows from (\ref{star}), (\ref{ws}) that $\mathrm{Pol}(Y)_q$ is a
$*$-subalgebra.

\medskip

\begin{lemma}\label{irry}
$\mathscr{L}$ is an irreducible $\mathrm{Pol}(Y)_q$-module.
\end{lemma}

\smallskip

{\bf Proof.} Let
$$
\mathscr{L}^{\mathrm{vac}}\stackrel{\mathrm{def}}{=}\bigcap_{J \ne
\{1,2,\dots,n
\}}\mathrm{Ker}\widetilde{\mathcal{T}}\left(t_{\{m+1,m+2,\dots,N
\}J}^{\wedge n}\right).
$$
$\mathscr{L}^{\mathrm{vac}}$ is invariant with respect to
$\widetilde{\mathcal{T}}(x)$ due to the commutation relations (\ref{xtij}).
It follows that $\mathscr{L}^{\mathrm{vac}}= \bigoplus
\limits_j(\mathscr{L}^{\mathrm{vac}}\cap \mathscr{L}_j)$.\footnote{See
\cite[Proposition 1.5]{Kac}.} Hence by lemma \ref{c=1},
$\mathscr{L}^{\mathrm{vac}}=\mathscr{L}_0$. On the other hand, by
proposition \ref{*sa}, $\mathscr{L}_0=\mathbb{C}\mathbf{e}_0$. So, we
conclude that $\mathscr{L}^{\mathrm{vac}}=\mathbb{C}\mathbf{e}_0$.

Turn back to the proof of irreducibility. Assume that there exists a
non-trivial $\mathrm{Pol}(Y)_q$-invariant subspace $\mathscr{L}'$. Then
$\mathscr{L}''=\bigoplus \limits_j(\mathscr{L}_j \ominus(\mathscr{L}'\cap
\mathscr{L}_j ))$ is also an invariant subspace and
$\mathscr{L}=\mathscr{L}'\oplus \mathscr{L}''$. Now apply the operators
$\widetilde{\mathcal{T}}\left(t_{\{m+1,m+2,\dots,N \}J}^{\wedge n}\right)$,
$J \ne \{1,2,\dots,n \}$ to find non-zero vectors from $\mathscr{L}'\cap
\mathscr{L}^{\mathrm{vac}}$, $\mathscr{L}''\cap \mathscr{L}^{\mathrm{vac}}$,
which are linear independent. This contradicts $\dim
\mathscr{L}^{\mathrm{vac}}=1$. \hfill $\square$

\smallskip

Turn back to proving that $\mathscr{J}$ is one-to-one. It follows from
lemmas \ref{Izm} and \ref{irry} that $\mathscr{J}$ is onto. To see that
$\mathscr{J}$ is injective, observe that $\mathscr{J}\mathcal{H}_i \subset
\mathscr{L}_i$, $\dim
\mathbb{C}[\mathrm{Mat}_{m,n}]_{q,i}=\dbinom{mn+i-1}{i}$, $\dim
\mathscr{L}_i=\dbinom{mn+i-1}{i}$, $i \in \mathbb{Z}_+$ (the latter equality
is due to proposition \ref{*sa}). \hfill $\square$

\medskip

{\sc Remark.} It follows from proposition \ref{iso} and lemmas \ref{Izm},
\ref{irry} that $T$ is irreducible.

\bigskip

\section{Boundedness of the quantum matrix ball}\label{bouT}

We use here the norm of an $m \times n$ matrix with entries in
$\mathrm{End}\,\mathcal{L}$ defined as the norm of the associated linear map
$\bigoplus \limits_{a=1}^n \mathcal{L}\to \bigoplus \limits_{\alpha=1}^m
\mathcal{L}$. Consider the matrices $\mathbf{Z}=(z_{\alpha
a})_{\alpha=1,\dots,m,a=1,\dots,n}$ with $z_{\alpha
a}=(-q)^{\alpha-1}z_a^{m+1-\alpha}$, and $\mathcal{T}(\mathbf{Z})=
(\mathcal{T}(z_{\alpha a}))_{\alpha=1,\dots,m,\;a=1,\dots,n}$.

\smallskip

\begin{proposition}\label{bou1}
$\|\mathcal{T}(\mathbf{Z})\|\le 1$.
\end{proposition}

We need the following

\begin{lemma}
In the matrix algebra with entries from $\mathbb{C}[SL_N]_{q,t}$ one has
\begin{equation}\label{IZ}
i(\mathbf{Z})=\mathbf{T}_{12}^{-1}\mathbf{T}_{11},
\end{equation}
with $i(\mathbf{Z})=(i(z_{\alpha a}))$, $\mathbf{T}_{11}=(t_{\alpha a})$,
$\mathbf{T}_{12}=(t_{\alpha,n+\beta})$, $\alpha,\beta=1,\dots,m$,
$a=1,2,\dots,n$.
\end{lemma}

\smallskip

{\bf Proof.} Let $\mathbf{t}=(t_{ij})_{i,j=1,\dots,m}$ and
$\det_q'\mathbf{t}=\sum_{s \in S_m}(-q)^{-l(s)}t_{s(m)m}t_{s(m-1)m-1}\cdots
t_{s(1)1}$. It is well known \cite[section 4]{PW} that in
$\mathbb{C}[\mathrm{Mat}_{m,m}]_q$ one has
\begin{equation}\label{det'}
\det \nolimits_q'\mathbf{t}=\det \nolimits_q \mathbf{t}.
\end{equation}
Now (\ref{IZ}) is derivable from (\ref{det'}) and the explicit form of
$\mathbf{T}^{-1}$ in the matrix algebra with entries from
$\mathbb{C}[\mathrm{Mat}_{m,m}]_q$ (see \cite{NYM}):
\begin{equation}\label{im}
\det \nolimits_q \mathbf{T}\cdot(\mathbf{T}^{-1})_{\alpha \beta}
=(-q)^{\alpha-\beta}\det \nolimits_q(\mathbf{T}_{\beta \alpha}),\qquad
\alpha,\beta=1,\dots,m
\end{equation}
(Here, just as in the classical case $q=1$, $\mathbf{T}_{\beta \alpha}$ is a
matrix derived from $\mathbf{T}$ by discarding the row $\beta$ and column
$\alpha$.) \hfill $\square$

\smallskip

{\bf Proof} of proposition \ref{bou1}. Apply (\ref{im}) to invert the matrix
$\mathbf{t}=(t_{ij})_{i,j=1,\dots,N}$ in the matrix algebra with entries
from $\mathbb{C}[SL_N]_{q,t}$:
$$
\sum_{a=1}^N(-q)^{a-\beta}t_{\alpha a}\det \nolimits_q(\mathbf{t}_{\beta
a})=\delta_{\alpha \beta},\qquad \alpha,\beta=1,\dots,N.
$$
Hence
$$
-\sum_{c=1}^nt_{\alpha c}t_{\beta
c}^*+\sum_{\gamma=1}^mt_{\alpha,n+\gamma}t_{\beta,n+\gamma}^*=\delta_{\alpha
\beta},\qquad \alpha,\beta=1,\dots,m.
$$
After introducing the notation
\begin{align*}
\mathbf{T}_{11}&=(t_{\alpha a})_{\alpha=1,\dots,m,\;a=1,\dots,n},&
\mathbf{T}_{12}&=(t_{\alpha,n+\beta})_{\alpha,\beta=1,\dots,m},
\\ \mathbf{T}_{11}^*&=(t_{a \alpha}^*)_{\alpha=1,\dots,m,\,a=1,\dots,n},&
\mathbf{T}_{12}^*&=(t_{n+\beta,\alpha}^*)_{\alpha,\beta=1,\dots,m},
\end{align*}
we get
\begin{equation}\label{bftt*}
-\mathbf{T}_{11}\mathbf{T}_{11}^*+\mathbf{T}_{12}\mathbf{T}_{12}^*=I.
\end{equation}
It follows from (\ref{bftt*}) and (\ref{IZ}) that
$i(I-\mathbf{Z}\mathbf{Z}^*)=\mathbf{T}_{12}^{-1}(\mathbf{T}_{12}^{-1})^*$.
Apply the representation $\widetilde{\mathcal{T}}$ to both parts of the
above relation. By a virtue of $\mathcal{T}=\widetilde{\mathcal{T}}\circ i$
we obtain $\mathcal{T}(I-\mathbf{Z}\mathbf{Z}^*)=
\widetilde{\mathcal{T}}(\mathbf{T}_{12}^{-1})
\widetilde{\mathcal{T}}(\mathbf{T}_{12}^{-1})^*\ge 0$. Hence
$\mathcal{T}(\mathbf{Z})\mathcal{T}(\mathbf{Z})^*\le I$,
$\|\mathcal{T}(\mathbf{Z})\|=\|\mathcal{T}(\mathbf{Z}^*)\|\le 1$. \hfill
$\square$

\medskip

Proposition \ref{bou1} and the unitary equivalence of the representations
$T$ and $\mathcal{T}$ (see section \ref{ue}) imply

\smallskip

\begin{corollary}\label{bou}
The operators ${T}(z_a^\alpha)$, $a=1,\dots,n$, $\alpha=1,\dots,m$, are
bounded.
\end{corollary}

\bigskip

\section{The quantum universal enveloping algebra \boldmath $U_q
\mathfrak{sl}_N$}\label{PBW}

The Drinfeld-Jimbo quantum universal enveloping algebra is among the basic
notions of the quantum group theory. Recall the definition of the Hopf
algebra $U_q \mathfrak{sl}_N$ \cite{Jant}. Let $(a_{ij})_{i,j=1,\ldots,N-1}$
be the Cartan matrix of $\mathfrak{sl}_N$:
\begin{equation}\label{Cm}
a_{ij}=
\begin{cases}
2,& i-j=0
\\ -1,& |i-j|=1
\\ 0,& \mathrm{otherwise}.
\end{cases}
\end{equation}
The algebra $U_q \mathfrak{sl}_N$ is determined by the generators $E_i$,
$F_i$, $K_i$, $K_i^{-1}$, $i=1,\ldots,N-1$, and the relations
\begin{gather}
K_iK_j=K_jK_i,\quad K_iK_i^{-1}=K_i^{-1}K_i=1,\quad
K_iE_j=q^{a_{ij}}E_jK_i,\quad K_iF_j=q^{-a_{ij}}F_jK_i \nonumber
\\ E_iF_j-F_jE_i=\delta_{ij}\,(K_i-K_i^{-1})/(q-q^{-1})\nonumber
\\ E_i^2E_j-(q+q^{-1})E_iE_jE_i+E_jE_i^2=0,\qquad |i-j|=1
\\ F_i^2F_j-(q+q^{-1})F_iF_jF_i+F_jF_i^2=0,\qquad |i-j|=1 \nonumber
\\ [E_i,E_j]=[F_i,F_j]=0,\qquad |i-j|\ne 1.\nonumber
\end{gather}
The comultiplication $\Delta$, the antipode $S$, and the counit
$\varepsilon$ are determined by
\begin{gather}
\Delta(E_i)=E_i \otimes 1+K_i \otimes E_i,\quad \Delta(F_i)=F_i \otimes
K_i^{-1}+1 \otimes F_i,\quad \Delta(K_i)=K_i \otimes K_i,\label{comult}
\\ S(E_i)=-K_i^{-1}E_i,\qquad S(F_i)=-F_iK_i,\qquad S(K_i)=K_i^{-1},
\\ \varepsilon(E_i)=\varepsilon(F_i)=0,\qquad \varepsilon(K_i)=1.\nonumber
\end{gather}

We consider in the sequel only $U_q \mathfrak{sl}_N$-modules of the form
$V=\bigoplus \limits_{\boldsymbol{\mu} \in
\mathbb{Z}^{N-1}}V_{\boldsymbol{\mu}}$, with
$\boldsymbol{\mu}=(\mu_1,\ldots,\mu_{N-1})$, $V_ {\boldsymbol{\mu}}=\{v \in
V|\,K_iv=q^{\mu_i}v,\,i=1,\ldots,N-1 \}$, to be referred to as weight
modules.\footnote{Note that some authors use the term 'weight' for a larger
class of modules \cite{KlSh}.} This agreement allows one to introduce the
linear operators in $V$
$$
H_jv=\mu_jv,\quad v \in V_{\boldsymbol{\mu}}, \qquad j=1,\ldots,N-1.
$$

Note that the defining relations in the classical universal enveloping
algebra $U \mathfrak{sl}_N$ can be derived from those in $U_q
\mathfrak{sl}_N$ via the substitution $K_i^{\pm 1}=q^{\pm H_i}$ and the
formal passage to a limit as $q \to 1$ (e.g.
$\lim\limits_{q\to1}\frac{K_i-K_i^{-1}}{q-q^{-1}}=H_i$, $i=1,\ldots,N-1$).

Equip all the weight $U_q \mathfrak{sl}_N$-modules with the gradation $\deg
v=j \:\Leftrightarrow \: H_0v=2jv$, \index{gradation for $U_q
\mathfrak{sl}_N$-modules} where $H_0$ is the unique element of the standard
Cartan subalgebra of $\mathfrak{sl}_N$ with the following properties:
$$[H_0,E_j]=0, \quad j \ne n; \qquad [H_0,E_n]=2E_n.$$
Let us present an explicit formula for $H_0$:
\begin{equation}\label{H0}
H_0=\frac2{m+n}\left(m \sum_{j=1}^{n-1}jH_j+n \sum_{j=1}^{m-1}jH_{N-j}+mnH_n
\right).
\end{equation}
It is easy to prove that

\begin{lemma}\label{orth}
$H_0$ is orthogonal to all the vectors $H_j$, $j \ne n$, with respect to an
invariant bilinear form in $\mathfrak{sl}_N$.
\end{lemma}

\smallskip

The rest of this section is intended to recall some well known results of
quantum group theory \cite{Jant}.

Recall that for the standard system of simple roots $\{\alpha_i
\}_{i=1,\ldots,N-1}$, of $\mathfrak{sl}_N$ one has $\alpha_i(H_j)=a_{j\,i}$,
$i,j=1,\ldots,N-1$, with $(a_{ij})$ being the Cartan matrix (\ref{Cm}). The
Weyl group is generated by simple reflections
$s_i(\alpha_j)=\alpha_j-a_{ij}\alpha_i$. In our case it is canonically
isomorphic to $S_N$: $s_i \mapsto(i,i+1)$. Consider the longest element
$w_0=(N,N-1,\ldots,2,1)\in S_N$, together with its reduced expression
$w_0=s_{i_1}\cdot s_{i_2}\cdot \ldots \cdot s_{i_M}$, $M=N(N-1)/2$, $1 \le
i_k \le N-1$. One can associate to the reduced expression a total order on
the set of positive roots of $\mathfrak{sl}_N$, and then a basis in the
vector space $U_q \mathfrak{sl}_N$. The total order is given by
$$
\beta_1=\alpha_1,\quad \beta_2=s_{i_1}(\alpha_{i_2}),\quad
\beta_3=s_{i_1}s_{i_2}(\alpha_{i_3}),\quad \ldots \quad\beta_M=s_{i_1}\ldots
s_{i_{M-1}}(\alpha_{i_M}).
$$

Turn to description of the basis in $U_q \mathfrak{sl}_N$ associated to the
reduced expression of $w_0$. G. Lusztig \cite{Lus} has defined an action of
the braid group $B_N$ as a group of automorphisms of the algebra $U_q
\mathfrak{sl}_N$ (we follow the definition given in \cite{DamDCon}):
\begin{align*}\label{la}
T_i(E_i)&=-F_iK_i, & T_i(F_i)&=-K_i^{-1}E_i,
\\ T_i(E_j)&=
\begin{cases}
E_j,& |i-j|>1,
\\ q^{-1}E_jE_i-E_iE_j,& |i-j|=1,
\end{cases}
& T_i(F_j)&=
\begin{cases}
F_j,& |i-j|>1,
\\ qF_iF_j-F_jF_i,& |i-j|=1,
\end{cases}
\\ T_i(K_j)&=K_jK_i^{-a_{ij}}.
\end{align*}
Note that the automorphisms $T_i$ permute the {\sl weight} spaces
$$
(U_q \mathfrak{g})_{\boldsymbol{\lambda}}=\{\xi \in U_q \mathfrak{g}|\,K_i
\xi K_i^{-1}=q^{\lambda_i}\xi \},\qquad
\boldsymbol{\lambda}=(\lambda_1,\lambda_2,\ldots,\lambda_{N-1})\in
\mathbb{Z}^{N-1},
$$
in the following way:
\begin{equation}\label{tsi}
T_i:(U_q \mathfrak{g})_{\boldsymbol{\lambda}}\to (U_q
\mathfrak{g})_{s_i(\boldsymbol{\lambda})}.
\end{equation}
Furthermore,
$$
T_i(U_q \mathfrak{sl}_n \otimes U_q \mathfrak{sl}_m)=U_q \mathfrak{sl}_n
\otimes U_q \mathfrak{sl}_m,\qquad i \ne n,
$$
with $U_q \mathfrak{sl}_n \subset U_q \mathfrak{sl}_N$ being the Hopf
subalgebra generated by $E_i$, $F_i$, $K_i^{\pm 1}$, $i=1,2,\ldots,n-1$, and
$U_q \mathfrak{sl}_m \subset U_q \mathfrak{sl}_N$ the Hopf subalgebra
generated by $E_{n+i}$, $F_{n+i}$, $K_{n+i}^{\pm 1}$, $i=1,2,\ldots,m-1$.

We have two maps $\alpha_i \mapsto E_i$, $\alpha_i \mapsto F_i$
$i=1,\ldots,N-1$. These maps, defined on the set of simple roots, are
extended onto the set of all positive roots as follows:
$E_{\beta_s}=T_{i_1}T_{i_2}\ldots T_{i_{s-1}}(E_{i_s})$,
$F_{\beta_s}=T_{i_1}T_{i_2}\ldots T_{i_{s-1}}(F_{i_s})$. We use below the
notation $U_q \mathfrak{n}_+$ (respectively, $U_q \mathfrak{n}_-$) for the
subalgebra in $U_q \mathfrak{sl}_N$ generated by $\{E_i \}$ (respectively,
$\{F_i \}$), $i=1,2,\ldots,N-1$.

\medskip

\begin{proposition}\label{la}\hfill
\\ i) $E_{\beta_1}^{k_1}\cdot E_{\beta_2}^{k_2}\cdot \ldots \cdot
E_{\beta_M}^{k_M}$, $(k_1,k_2,\ldots,k_M)\in \mathbb{Z}_+^M$, constitute a
basis of weight vectors in the vector space $U_q \mathfrak{n}_+$;
\\ ii) $F_{\beta_M}^{j_M}\cdot F_{\beta_{M-1}}^{j_{M-1}}\cdot \ldots \cdot
F_{\beta_1}^{j_1}$, $(j_1,j_2,\ldots,j_M)\in \mathbb{Z}_+^M$, constitute a
basis of weight vectors in the vector space $U_q \mathfrak{n}_-$;
\\ iii) $F_{\beta_M}^{j_M}\cdot F_{\beta_{M-1}}^{j_{M-1}}\cdot \ldots \cdot
F_{\beta_1}^{j_1}\cdot K_1^{i_1}\cdot K_2^{i_2}\cdot \ldots \cdot
K_{N-1}^{i_{N-1}}\cdot E_{\beta_1}^{j_1}\cdot E_{\beta_2}^{j_2}\cdot \ldots
\cdot E_{\beta_M}^{j_M}$,
\\ $(k_1,k_2,\ldots,k_M),\;(j_1,j_2,\ldots,j_M)\in \mathbb{Z}_+^M$,
$(i_1,i_2,\ldots,i_{N-1})\in \mathbb{Z}^{N-1}$,
\\ constitute a basis of weight vectors in the vector space $U_q \mathfrak{sl}_N$.
\end{proposition}

\medskip

Consider a $U_q \mathfrak{sl}_N$-module $V^h$, determined by its generator
$v^h$ and the relations
\begin{gather*}
E_jv^h=(K_j^{\pm 1}-1)v^h=0,\qquad j=1,\ldots,N-1,
\\ F_iv^h=0,\qquad i=1,\ldots,n-1,n+1,\ldots,N-1.
\end{gather*}
($V^h$ is canonically isomorphic to the generalized Verma module with zero
highest weight.) We are about to apply proposition \ref{la} to produce a
basis of the vector space $V^h$ formed by homogeneous vectors. For that, we
use the class of reduced expressions for the element $w_0$ described below.
Recall the notation $M=N(N-1)/2$, $M'=M-mn$, $s_j=(j,j+1)$. Consider the
longest element for the subgroup $S_n \times S_m \subset S_N$
$$w_0'=(n,n-1,\ldots,1,N,N-1,\ldots,n+1),$$
together with the permutation
$$w''_0=(n+1,n+2,\ldots,N-1,N,1,2,\ldots,n-1,n).$$
Obviously, $w_0=w_0'\cdot w_0''$. Fix the reduced expression
$w_0=s_{i_1}s_{i_2}\ldots s_{i_M}$, given by concatenation of reduced
expressions for $w_0'$ and $w_0''$. It follows from the definitions that
$\deg \left(F_{\beta_j}\right)=-1$ for $j>M'$, and $\deg
\left(F_{\beta_j}\right)=0$ for $j \le M'$. Thus the vectors
$$
F_{\beta_M}^{k_M}F_{\beta_{M-1}}^{k_{M-1}} \ldots
F_{\beta_{M'+1}}^{k_{M'+1}}\cdot v^h,\qquad (k_{M'+1},\ldots,k_M)\in
\mathbb{Z}_+^{mn},
$$
constitute a basis in the vector space $V^h$, and
$$
\deg \left(F_{\beta_M}^{k_M}F_{\beta_{M-1}}^{k_{M-1}} \ldots
F_{\beta_{M'+1}}^{k_{M'+1}}\cdot v^h \right)=-\sum_{j=M'+1}^{M}k_j.
$$

It is easy to obtain a similar result for a weight $U_q
\mathfrak{sl}_N$-module $V^l$ determined by its generator $v^l$ and the
relations
\begin{gather*}
F_jv^l=(K_j^{\pm 1}-1)v^l=0,\qquad j=1,\ldots,N-1,
\\ E_iv^l=0,\qquad i=1,\ldots,n-1,n+1,\ldots,N-1.
\end{gather*}
The following vectors form a basis of the graded vector space $V^l$
consisting of homogeneous vectors:
$$
S \left(E_{\beta_{M'+1}}^{j_{M'+1}}\cdot E_{\beta_{M'+2}}^{j_{M'+2}}\cdot
\ldots \cdot E_{\beta_M}^{j_M}\right)v^l,
$$
with $(j_{M'+1},j_{M'+2},\ldots,j_M)\in \mathbb{Z}_+^{mn}$, and $S$ being
the antipode of the Hopf algebra $U_q \mathfrak{sl}_N$.

Recall some notions of the theory of Hopf algebras \cite{CP}. Let $A$ be an
abstract Hopf algebra and $F$ an algebra equipped also with a structure of
$A$-module. $F$ is said to be an $A$-module algebra if the multiplication $F
\otimes F \to F$, $f_1 \otimes f_2 \mapsto f_1f_2$, is a morphism of
$A$-modules. In the case of a unital algebra $F$, the additional assumption
is introduced that the embedding $\mathbb{C}\hookrightarrow F$, $1 \mapsto
1$, is a morphism of $A$-modules. A duality argument allows one also to
introduce a notion of $A $-module coalgebra.


The results, cited below, are due to S. Levendorskii and Ya. Soibelman
\cite{LS1, LS2} and A. Kirillov, N. Reshetikhin \cite{KR, CP}. For a good
survey of those the reader is referred to \cite{R, DamDCon}.

Let $V_1$ and $V_2$ be some $U_q\mathfrak{sl}_{N}$-modules. It is well known
that, in general, the ordinary flip $$\sigma_{V_1, V_2}: V_1\otimes V_2\to
V_2\otimes V_1, \quad \sigma_{V_1, V_2}: v_1\otimes v_2\to v_2\otimes v_1$$
is not a morphism of $U_q\mathfrak{sl}_{N}$-modules. V. Drinfeld \cite{Dr1}
introduced the extremely important notion of {\it the universal $R$-matrix}
which has lead to appropriate $q$-analogs of the operators $\sigma_{V_1,
V_2}$. Let us describe these $q$-analogues.

To start with, recall the standard notation $U_q\mathfrak{b}^+$
(respectively $U_q\mathfrak{b}^-$) for the Hopf subalgebra in
$U_q\mathfrak{sl}_N$ generated by $K_i^{\pm1}$, $E_i$, $i=1,\ldots,N-1$
(respectively $K_i^{\pm1}$, $F_i$, $i=1,\ldots,N-1$). We denote by
$\mathcal{C}^+$ (respectively $\mathcal{C}^-$) the category of
$U_q\mathfrak{b}^+$-locally finite dimensional (respectively
$U_q\mathfrak{b}^-$-locally finite dimensional) weight
$U_q\mathfrak{sl}_{N}$-modules.

Let $V_1$, $V_2$ be weight $U_q \mathfrak{sl}_N$-modules, and either $V_1\in
\mathcal{C}^+$ or $V_2\in \mathcal{C}^-$. 
The formula below determines a linear operator $R_{V_1,V_2}$ in $V_1 \otimes
V_2$:
\begin{multline}\label{rmmf}
R=\exp_{q^2}\left((q^{-1}-q)E_{\beta_M}\otimes F_{\beta_M}\right)\cdot
\exp_{q^2}\left((q^{-1}-q)E_{\beta_{M-1}}\otimes F_{\beta_{M-1}}\right)\cdot
\ldots \cdot
\\ \cdot \exp_{q^2}\left((q^{-1}-q)E_{\beta_1}\otimes
F_{\beta_1}\right)q^{-t_0},
\end{multline}
with $\exp_{q^2}(u)=\sum \limits_{k=0}^\infty \dfrac{u^k}{(k)_{q^2}!}$;
$(k)_{q^2}!=\prod \limits_{j=1}^k\dfrac{1-q^{2j}}{1-q^2}$,
\begin{equation}\label{t0}
t_0=\sum \limits_{i,j=1}^{N-1}c_{ij}H_i \otimes H_j,
\end{equation}
and $(c_{ij})_{i,j=1,\ldots,N-1}$ is inverse to the Cartan matrix
$(a_{ij})_{i,j=1,\ldots,N-1}$. It is worthwhile to note that $c_{ij}\in
\dfrac1N \mathbb{Z}_+$.

Now we use the relation $\alpha_i(H_j)=a_{j\,i}$, $i,j=1,\ldots,N-1$, to get
a different description of $t_0$:
$$\alpha_i \otimes \alpha_j(t_0)=a_{j\,i},\qquad i,j=1,\ldots,N-1.$$
Recall the definition of the standard inner product in the Cartan
subalgebra: $(H_i,H_j)=a_{ij}$, $i,j=1,\ldots,N-1$. It allows one to get the
third description of $t_0$:
$$(t_0,H_i \otimes H_j)=(H_i,H_j);\qquad i,j=1,\ldots,N-1.$$
That is,
\begin{equation}\label{3descr}
t_0=\sum \limits_{k=1}^{N-1}\dfrac{I_k \otimes I_k}{(I_k,I_k)}
\end{equation}
for any orthogonal basis of the Cartan subalgebra.

The formula (\ref{rmmf}) involves analogs of root vectors of the Lie algebra
$\mathfrak{sl}_N$ whose construction, we recall, depends on the choice of a
reduced expression for the longest element $w_0\in S_N$. Nevertheless, it is
well known that the operators $R_{V_1,V_2}$ are independent of that choice.

 Let us list some properties of the operators
$\check{R}_{V_1,V_2}\stackrel{\rm def}{=}\sigma_{V_1,V_2}\cdot R_{V_1,V_2}$.
Let again $V_1$, $V_2$ be weight $U_q \mathfrak{sl}_N$-modules, and either
$V_1\in \mathcal{C}^+$ or $V_2\in \mathcal{C}^-$. Then the operator
$\check{R}_{V_1,V_2}:V_1\otimes V_2\to V_2\otimes V_1$ is an invertible
operator and a morphism of $U_q \mathfrak{sl}_N$-modules.


Suppose $V$, $V_1$, $V_2$ are weight $U_q \mathfrak{sl}_N$-modules, and
either $V_1, V_2\in \mathcal{C}^+$ or $V\in \mathcal{C}^-$. Then
\begin{equation}\label{braid1}
\check{R}_{V_1\otimes V_2,
V}=(\check{R}_{V_1,V}\otimes\mathrm{id}_{V_2})\cdot(\mathrm{id}_{V_1}\otimes
\check{R}_{V_2,V}).
\end{equation}

Finally, suppose $V$, $V_1$, $V_2$ are weight $U_q \mathfrak{sl}_N$-modules,
and either $V\in \mathcal{C}^+$ or $V_1, V_2\in \mathcal{C}^-$. Then
\begin{equation}\label{braid2}
\check{R}_{V,V_1\otimes V_2}=(\mathrm{id}_{V_1}\otimes
\check{R}_{V,V_2})\cdot(\check{R}_{V,V_1}\otimes\mathrm{id}_{V_2}).
\end{equation}

The above properties of the operators $\check{R}_{V_1, V_2}$ allow one to
treat them as $q$-analogues of the ordinary flips $\sigma_{V_1,V_2}$.

Consider the vector representation $\pi$ of $U_q \mathfrak{sl}_N$ in
$\mathbb{C}^N$:
\begin{align*}
\pi(E_i)e_j&=
\begin{cases}
q^{-1/2}e_{j-1},& j=i+1
\\ 0,& \mathrm{otherwise}
\end{cases}
& \pi(F_i)e_j&=
\begin{cases}
q^{1/2}e_{j+1},& j=i
\\ 0,& \mathrm{otherwise}
\end{cases}
\\ \pi(K_i^{\pm 1})e_j&=
\begin{cases}
q^{\pm 1}e_j,& j=i
\\ q^{\mp 1}e_j,& j=i+1
\\ e_j,& \mathrm{otherwise}
\end{cases}
\end{align*}
with $i=1,2,\ldots,N-1$, $j=1,2,\ldots,N$, $\{e_j\}$ being the standard
basis in $\mathbb{C}^N$. The linear functionals $l_{j\,k}\in (U_q
\mathfrak{sl}_N)^*$ given by
$$
\pi(\xi)e_k=\sum_{j=1}^Nl_{jk}(\xi)e_j,\qquad \xi \in U_q \mathfrak{sl}_N,
$$
are called matrix elements of $\pi$ with respect to the basis $\{e_j\}$.
There exists a canonical non-degenerate pairing (see e.g. \cite{VaSo2})
$$
\mathbb{C}[SL_N]_q \times U_q \mathfrak{sl}_N \to \mathbb{C},\qquad f \times
\xi \mapsto \langle f,\xi \rangle,
$$
which determines an embedding of the Hopf algebras
\begin{equation}\label{local1}
\mathbb{C}[SL_N]_q \hookrightarrow (U_q \mathfrak{sl}_N)^*,\qquad
t_{jk}\mapsto l_{jk},\quad j,k=1,\ldots,N.
\end{equation}

Let $L(\lambda)$ be the simple finite dimensional weight $U_q
\mathfrak{sl}_N$-module with highest weight $\lambda$. The embedding
$(\mathrm{End}_\mathbb{C}L(\lambda))^*\hookrightarrow(U_q
\mathfrak{sl}_N)^*$ allows one to get an isomorphism
$$
\mathbb{C}[SL_N]_q \simeq \bigoplus_\lambda
(\mathrm{End}_\mathbb{C}L(\lambda))^*.
$$
The embedding (\ref{local1}) may be used to equip $\mathbb{C}[SL_N]_q$ with
a structure of $U_q \mathfrak{sl}_N$-module algebra:
$$
\langle \xi f,\eta \rangle=\langle f,\eta \xi \rangle,\qquad f \in
\mathbb{C}[SL_N]_q,\quad \xi,\eta \in U_q \mathfrak{sl}_N.
$$
It is now deducible from the definitions that the generators of $U_q
\mathfrak{sl}_N$ act on the generators of $\mathbb{C}[SL_N]_q$ in the
following way:
\begin{gather}
E_it_{j,k}=
\begin{cases}
q^{-1/2}t_{j,k-1},& k=i+1
\\ 0,& \mathrm{otherwise}
\end{cases},\qquad
F_it_{j,k}=
\begin{cases}
q^{1/2}t_{j,k+1},& k=i
\\ 0,& \mathrm{otherwise}
\end{cases},\label{EFt}
\\ K_i^{\pm 1}t_{j,k}=
\begin{cases}
q^{\pm 1}t_{j,k},& k=i
\\ q^{\mp 1}t_{j,k},& k=i+1
\\ t_{j,k},& \mathrm{otherwise}.
\end{cases}\label{Kt}
\end{gather}

\bigskip

\section{The algebras \boldmath $\mathbb{C}[\mathcal{M}at_{m,n}]_q$,
$\mathrm{Pol}(\mathcal{M}at_{m,n})_q$}\label{da}

In all the above observations we assumed that $q \in(0,1)$ and the ground
field is $\mathbb{C}$. Nevertheless, it appears to be much more convenient
in this section to replace $\mathbb{C}$ with the field $\mathbb{C}(q^{1/s})$
of rational functions of the indeterminate $q^{1/s}$, $s=2N$.\footnote{One
can observe from the formulation of proposition \ref{zaalph} that $s$ should
be even, and $s \in N \mathbb{Z}$ due to (\ref{rhat}), (\ref{rmmf}),
(\ref{t0})} In the subsequent sections we are going to retrieve our original
convention concerning the ground field.

Since our goals are results with $\mathbb{C}$ as a ground field, we need an
appropriate procedure for backward passage from $\mathbb{C}(q^{1/s})$ to a
ring of Laurent polynomials and finally to $\mathbb{C}$. A passage of that
kind could be done via standard techniques\footnote{non-restricted
specialization} well known in quantum group theory \cite[\S 9.2]{CP}. In
what follows we obtain a number of results for algebras and $*$-algebras
over $\mathbb{C}(q^{1/s})$ (it is implicit in this context that
$\overline{q^{1/s}}=q^{1/s}$). We keep the former notations
$\mathbb{C}[\mathrm{Mat}_{m,n}]_q$, $\mathrm{Pol}(\mathrm{Mat}_{m,n})_q$,
$\mathbb{C}[SL_N]_q$, $\mathbb{C}[\widetilde{G}]_q$, $U_q \mathfrak{sl}_N$
for algebras over $\mathbb{C}(q^{1/s})$ determined by 'the same' generators
and relations as before in the case of the ground field $\mathbb{C}$. It is
well known that the results of section \ref{PBW} are valid also in the case
of the ground field $\mathbb{C}(q^{1/s})$ (cf. \cite{DamDCon}).

In this section we are going to develop a different approach to the algebra
$\mathbb{C}[\mathrm{Mat}_{m,n}]_q$. More precisely, we are about to
construct an algebra $\mathbb{C}[\mathcal{M}at_{m,n}]_q$ which is
canonically isomorphic to $\mathbb{C}[\mathrm{Mat}_{m,n}]_q$ and is much
more convenient for our further goals.

Consider the Hopf algebra $U_q \mathfrak{sl}_N^{\mathrm{op}}$ which differs
from $U_q \mathfrak{sl}_N$ by replacing its comultiplication by the opposite
one. Equip $V^h$ with a structure of $U_q
\mathfrak{sl}_N^{\mathrm{op}}$-module coalgebra: $\Delta:v^h \mapsto v^h
\otimes v^h$. Consider the graded vector space dual to $V^h$:
$$
\mathbb{C}[\mathcal{M}at_{m,n}]_q=\bigoplus_{j=0}^\infty
\mathbb{C}[\mathcal{M}at_{m,n}]_{q,j},\qquad
\mathbb{C}[\mathcal{M}at_{m,n}]_{q,j}=\left(V_{-j}^h\right)^*,\quad j \in
\mathbb{Z}_+.
$$
Equip $\mathbb{C}[\mathcal{M}at_{m,n}]_q$ with a structure of $U_q
\mathfrak{sl}_N$-module algebra by the duality:
$$
\langle \xi f,v \rangle=\langle f,S(\xi)v \rangle, \qquad \langle f_1f_2,v
\rangle=\sum \limits_i \langle f_1,v'_i\rangle \langle f_2,v''_i \rangle,
$$
with $f,f_1,f_2 \in \mathbb{C}[\mathcal{M}at_{m,n}]_q$, $v \in V^h$,
$\triangle v=\sum \limits_jv'_j \otimes v''_j$.

Our immediate intention is to describe $\mathbb{C}[\mathcal{M}at_{m,n}]_q$
in terms of generators and relations. Consider the Hopf subalgebra $U_q
\mathfrak{sl}_n \subset U_q \mathfrak{sl}_N$ generated by $E_i$, $F_i$,
$K_i^{\pm 1}$, $i=1,2,\ldots,n-1$, and the Hopf subalgebra $U_q
\mathfrak{sl}_m \subset U_q \mathfrak{sl}_N$ generated by $E_{n+i}$,
$F_{n+i}$, $K_{n+i}^{\pm 1}$, $i=1,2,\ldots,m-1$. It follows from the
definitions that the homogeneous component
$\mathbb{C}[\mathcal{M}at_{m,n}]_{q,1}=\{f \in
\mathbb{C}[\mathcal{M}at_{m,n}]_q|\:\deg f=1 \}=(V^h_{-1})^*$ is a $U_q
\mathfrak{sl}_n \otimes U_q \mathfrak{sl}_m$-module. We are going to prove
that this module splits into the tensor product of a $U_q
\mathfrak{sl}_n$-module related to the vector representation and a $U_q
\mathfrak{sl}_m$-module related to the covector representation. Consider the
$U_q \mathfrak{sl}_n$-module $U$ and the $U_q \mathfrak{sl}_m$-module $V$,
determined in the bases $\{u_a \}_{a=1,\ldots,n}$, $\{v^\alpha
\}_{\alpha=1,\ldots,m}$ by
\begin{align*}
E_iu_a&=
\begin{cases}
q^{-1/2}u_{a-1},& a=i+1
\\ 0,& \mathrm{otherwise}
\end{cases}
&E_{n+i}v^\alpha&=
\begin{cases}
q^{-1/2}v^{\alpha-1},& \alpha=m-i+1
\\ 0,& \mathrm{otherwise}
\end{cases}
\\ F_iu_a&=
\begin{cases}
q^{1/2}u_{a+1},& a=i
\\ 0,& \mathrm{otherwise}
\end{cases}
& F_{n+i}v^\alpha&=
\begin{cases}
q^{1/2}v^{\alpha+1},& \alpha=m-i
\\ 0,& \mathrm{otherwise}
\end{cases}
\\ K_i^{\pm 1}u_a&=
\begin{cases}
q^{\pm 1}u_a,& a=i
\\ q^{\mp 1}u_a,& a=i+1
\\ u_a,& \mathrm{otherwise}
\end{cases}
& K_{n+i}^{\pm 1}v^\alpha&=
\begin{cases}
q^{\pm 1}v^\alpha,& \alpha=m-i
\\ q^{\mp 1}v^\alpha,& \alpha=m-i+1
\\ v^\alpha,& \mathrm{otherwise}
\end{cases}.
\end{align*}

\medskip

\begin{proposition}\label{zaalph}
There exists a unique collection $\{z_a^\alpha
\}_{a=1,\ldots,n;\:\alpha=1,\ldots,m}$, of elements of
$\mathbb{C}[\mathcal{M}at_{m,n}]_{q,1}$ such that the map $i:u_a \otimes
v^\alpha \mapsto z_a^\alpha$, $a=1,\ldots,n;\:\alpha=1,\ldots,m$ admits an
extension up to an isomorphism of $U_q \mathfrak{sl}_n \otimes U_q
\mathfrak{sl}_m$-modules $i:U \otimes V \mapsto
\mathbb{C}[\mathcal{M}at_{m,n}]_{q,1}$, and $F_nz_n^m=q^{1/2}$.
\end{proposition}

\smallskip

{\bf Proof.} Let $V^h_{-k}$ denotes the $(-k)$-th graded component of the
$U_q\mathfrak{sl}_N$-module $V^h$: $$V^h_{-k}{=}\{v|\:H_0v=-2kv \}.$$ It
follows from the results of the previous section that the elements
$F_{\beta_M}^{k_M}F_{\beta_{M-1}}^{k_{M-1}} \ldots
F_{\beta_{M'+1}}^{k_{M'+1}}v^h $, $k_{M'+1}+k_{M'+2}+\ldots+k_M=k$,
constitute a basis in $V^h_{-k}$.
 Hence, the dimension
of $V^h_{-k}$ is just the same as in the classical $(q=1)$ case:
\begin{equation}\label{dim}
\dim V^h_{-k}=\binom{mn+k-1}k.
\end{equation}
Observe that $V^h_{-1}$ is non-zero, so $v'=F_nv^h \ne 0$ and
$$
E_jv'=0,\qquad H_jv'=
\begin{cases}
-2v',& j=n
\\ v',& |j-n|=1
\\ 0,& |j-n|>1
\end{cases},
\qquad j=1,\ldots,N-1.
$$
This, together with $\dim V^h_{-1}=mn$ implies that $v'$ generates the
simple weight $U_q \mathfrak{sl}_n \otimes U_q \mathfrak{sl}_m$-module $
V^h_{-1}$, and $\mathbb{C}[\mathcal{M}at_{m,n}]_{q,1}\simeq U \otimes V $.
Of course, the isomorphism $i:U \otimes V \to
\mathbb{C}[\mathcal{M}at_{m,n}]_{q,1}$ is unique up to a multiple from the
ground field, and the elements $z_a^\alpha=i(u_a \otimes v^\alpha)$,
$a=1,\ldots,n$, $\alpha=1,\ldots,m$, satisfy all the requirements of our
proposition, except, possibly, the last property $F_nz_n^m=q^{1/2}$. One can
readily choose the above multiple in the definition of $i$, which provides
this property unless $F_nz_n^m=0$. In the latter case one has
$F_n(E_{i_1}^{k_1}E_{i_2}^{k_2}\ldots E_{i_l}^{k_l}z_n^m)=0$ for all
$i_1,\ldots,i_l$ different from $n$ and all $k_1, k_2,\ldots,k_l \in
\mathbb{Z}_+$. Hence it follows from the irreducibility of the $U_q
\mathfrak{sl}_n \otimes U_q \mathfrak{sl}_m$-module
$\mathbb{C}[\mathcal{M}at_{m,n}]_{q,1}\simeq U \otimes V$ that $F_n
\mathbb{C}[\mathcal{M}at_{m,n}]_{q,1}=0$, and thus $F_nv^h=0$. That is,
$\dim V^h=1$. On the other hand, it follows from (\ref{dim}) that $\dim
V^h=\infty$. This contradiction shows that $F_nz_n^m \ne 0$. \hfill
$\square$

\medskip

\begin{proposition}\label{gen}
$z_a^\alpha$, $a=1,\ldots,n$, $\alpha=1,\ldots,m$, generate the algebra
$\mathbb{C}[\mathcal{M}at_{m,n}]_q$.
\end{proposition}

\smallskip

{\bf Proof.} By a virtue of (\ref{dim}), it suffices to prove that for any
$k\in\mathbb{Z}_+$ the $\dbinom{mn+k-1}k$ monomials
\begin{equation}\label{mono}
(z_1^1)^{j_1^1}(z_1^2)^{j_1^2}\ldots(z_n^m)^{j_n^m},\qquad
j_1^1+j_1^2+\ldots+j_n^m=k
\end{equation}
 are linearly independent in $\mathbb{C}[\mathcal{M}at_{m,n}]_{q,k}$.
 An application of the standard techniques of specialization
(see \cite[chapter 5]{Jant}) allows us to reduce this statement to its
classical analogue. Namely, consider the basis
\begin{equation}\label{basisel}
F_{\beta_M}^{k_M}F_{\beta_{M-1}}^{k_{M-1}} \ldots
F_{\beta_{M'+1}}^{k_{M'+1}}v^h, \qquad k_M+k_{M-1}+\ldots+k_{M'}=k
\end{equation}
in $V^h_{-k}$. 
Let us denote by $\langle\,,\,\rangle$ the pairing
\begin{equation}\label{pairi}
\mathbb{C}[\mathcal{M}at_{m,n}]_{q}\times V^h\to\mathbb{C}(q^{1/s})
\end{equation}
 which is implicit in the definition of
$\mathbb{C}[\mathcal{M}at_{m,n}]_{q}$. Clearly, to prove linear independence
of the monomials (\ref{mono}), it suffices to show that the determinant of
the pairing
$\langle\,,\,\rangle:{\mathbb{C}[\mathcal{M}at_{m,n}]_{q,k}\times
V^h_{-k}}\to\mathbb{C}(q^{1/s})$
in the bases (\ref{mono}) and (\ref{basisel}) is non-zero.
 An idea, underlying the specialization techniques, may be described
roughly as follows. The determinant is a rational function of $q^{1/s}$. To
prove that the function is non-zero, it is enough to demonstrate that the
point $q=1$ is neither its pole nor its zero.

 Consider the ring
$\mathcal{A}=\mathbb{Q}[q^{1/s},q^{-1/s}]$ of Laurent polynomials in the
indeterminate $q^{1/s}$ and the $\mathcal{A}$-subalgebra $U_\mathcal{A}$ in
$U_q \mathfrak{sl}_N$ generated by the elements $E_i$, $F_i$, $K_i^{\pm 1}$,
$L_i=\dfrac{K_i-K_i^{-1}}{q-q^{-1}}$, $i=1,\ldots, N-1$. This is a Hopf
algebra:
$$\Delta(L_i)=L_i \otimes K_i+K_i^{-1}\otimes L_i, \quad S(L_i)=-L_i, \quad
\varepsilon(L_i)=0, \quad i=1,\ldots,N-1.$$ Let
$V_\mathcal{A}=U_\mathcal{A}v^h$. It is easy to show that the basis elements
$F_{\beta_M}^{k_M}F_{\beta_{M-1}}^{k_{M-1}} \ldots
F_{\beta_{M'+1}}^{k_{M'+1}}v^h$ are in $V_\mathcal{A}$. Denote by
$F_\mathcal{A} \subset \mathbb{C}[\mathcal{M}at_{m,n}]_{q}$ the
$\mathcal{A}$-module generated by all the monomials
$(z_1^1)^{j_1^1}(z_1^2)^{j_1^2}\ldots(z_n^m)^{j_n^m}$,
$j_1^1,j_1^2,\ldots,j_n^m\in\mathbb{Z}_+$. It follows from the relations in
$U_q \mathfrak{sl}_N$ and the definitions of modules $U$, $V$ that the value
of the linear functional $z_a^\alpha$ on a vector $v \in V_\mathcal{A}$ is
in $\mathcal{A}$. Hence, a similar statement is also valid for all $f \in
F_\mathcal{A}$. In particular, the aforementioned determinant belongs to
$\mathcal{A}$. We intend to prove that the determinant is a non-zero element
in $\mathcal{A}$. For that, it suffices to prove that its image under the
natural homomorphism $\mathcal{A}\to \mathbb{Q}$, $q^{1/s}\mapsto 1$ is a
non-zero number. The latter is a straightforward consequence of
non-degeneracy of the natural pairing
$$\mathbb{C}[z_1^1,\ldots,z_n^m]\times U\mathfrak{p}_-\to\mathbb{C}, \qquad
f(z_1^1,\ldots,z_n^m)\times\xi\mapsto S(\xi)(f(z_1^1,\ldots,z_n^m))|_{
z_1^1=\ldots=z_n^m=0},$$ where $U\mathfrak{p}_-$ is the universal enveloping
algebra of the Abelian Lie subalgebra
$$\mathfrak{p}_-=\{\xi\in\mathfrak{sl}_N\,|\,[H_0,\xi]=-2\xi\}$$ ($H_0$ is given
in (\ref{H0})), and $S$ is the antipode in $U\mathfrak{sl}_N$.\hfill
$\square$



\medskip

\begin{proposition}\label{zcr} The elements $z_a^\alpha$, $\alpha=1,\ldots, m$,
$a=1,\ldots, n$ of the algebra $\mathbb{C}[\mathcal{M}at_{m,n}]_q$ satisfy
the relations (\ref{zaa1}), (\ref{zaa2}), (\ref{zaa3}).
\end{proposition}

\smallskip

{\bf Proof.} It is easy to verify that the linear span of the left hand
sides in (\ref{zaa1}) -- (\ref{zaa3}) corresponds to a $U_q \mathfrak{sl}_n
\otimes U_q \mathfrak{sl}_m$-submodule of $M^{\otimes 2}$ with
$M=\mathbb{C}[\mathcal{M}at_{m,n}]_{q,1}$. Let $M_\mathcal{A} \subset M$ be
the $\mathcal{A}$-module generated by $\{z_a^\alpha \}$, $a=1,\ldots,n$,
$\alpha=1,\ldots,m$. By a virtue of proposition \ref{zaalph}, $M^{\otimes
2}$ is decomposed into a direct sum of four simple pairwise non-isomorphic
$U_q \mathfrak{sl}_n \otimes U_q \mathfrak{sl}_m$-modules.\footnote{The
tensor square of the vector (co-vector) representation is isomorphic to the
direct sum of its symmetric square and exterior square.} A similar
decomposition is also valid for $M_\mathcal{A} \otimes M_\mathcal{A}$, and a
specialization at $q=1$ leads to four pairwise non-isomorphic $U
\mathfrak{sl}_n \otimes U \mathfrak{sl}_m$-modules. By misuse of language,
one can say that each submodule of the $U_q \mathfrak{sl}_n \otimes U_q
\mathfrak{sl}_m$-module $M^{\otimes 2}$ is unambiguously determined by its
specialization at $q=1$. Consider two such submodules, namely, the kernel of
the multiplication operator $\mathbb{C}[\mathrm{Mat}_{m,n}]_{q,1}^{\otimes
2}\to \mathbb{C}[\mathrm{Mat}_{m,n}]_{q,2}$, $f_1 \otimes f_2 \mapsto
f_1f_2$ and the linear span of the elements given by the left hand sides of
(\ref{zaa1}) -- (\ref{zaa3}). Their specializations at $q=1$ coincide, and
hence the $U_q \mathfrak{sl}_n \otimes U_q \mathfrak{sl}_m$-submodules
themselves are the same. \hfill $\square$

\medskip

\begin{corollary}\label{hom}
There exists a unique homomorphism of graded algebras
\begin{equation}\label{zz}
j:\mathbb{C}[\mathrm{Mat}_{m,n}]_q\to
\mathbb{C}[\mathcal{M}at_{m,n}]_q,\qquad j:z_a^\alpha \mapsto
z_a^\alpha,\quad a=1,\ldots,n,\;\alpha=1,\ldots,m.
\end{equation}
\end{corollary}

\medskip

\begin{proposition}\label{isomo}
The homomorphism (\ref{zz}) is an isomorphism.
\end{proposition}

\smallskip

{\bf Proof.} It is an easy exercise to compute the dimensions of the
homogeneous components $\mathbb{C}[\mathrm{Mat}_{m,n}]_{q,k}=\{f \in
\mathbb{C}[\mathrm{Mat}_{m,n}]_q|\;\deg f=k \}$. Specifically,
\begin{equation}\label{dim1}
\dim \mathbb{C}[\mathrm{Mat}_{m,n}]_{q,k}=\binom{mn+k-1}k.
\end{equation}
It follows from proposition \ref{gen} that $j$ is onto. What remains is to
apply (\ref{dim}), (\ref{dim1}) to observe coincidence of the dimensions of
the graded components:
$$
\dim \mathbb{C}[\mathcal{M}at_{m,n}]_{q,k}=\dim
\mathbb{C}[\mathrm{Mat}_{m,n}]_{q,k},\qquad k \in \mathbb{Z}_+.\eqno \square
$$

\medskip

So far we considered the $U_q \mathfrak{sl}_N$-module algebra
$\mathbb{C}[\mathcal{M}at_{m,n}]_q$ dual to $V^h$. Now turn to producing a
$U_q \mathfrak{sl}_N$-module algebra
$\mathbb{C}[\overline{\mathcal{M}at}_{m,n}]_q$ dual to $V^l$. Equip $V^l$
with a structure of $U_q \mathfrak{sl}_N^{\mathrm{op}}$-module coalgebra:
$\Delta:v^l \mapsto v^l \otimes v^l$. Consider the graded vector space dual
to $V^l$:
$$
\mathbb{C}[\overline{\mathcal{M}at}_{m,n}]_q=\bigoplus_{j=0}^\infty
\mathbb{C}[\overline{\mathcal{M}at}_{m,n}]_{q,-j},\qquad
\mathbb{C}[\overline{\mathcal{M}at}_{m,n}]_{q,-j}=
\left(V_j^l\right)^*,\quad j \in \mathbb{Z}_+.
$$
Equip $\mathbb{C}[\overline{\mathcal{M}at}_{m,n}]_q$ with a structure of
$U_q \mathfrak{sl}_N$-module algebra by the duality:
$$
\langle \xi f,v \rangle=\langle f,S(\xi)v \rangle, \qquad \langle f_1f_2,v
\rangle=\sum \limits_i \langle f_1,v'_i\rangle \langle f_2,v''_i \rangle,
$$
with $f,f_1,f_2 \in \mathbb{C}[\overline{\mathcal{M}at}_{m,n}]_q$, $v \in
V^l$, $\triangle v=\sum \limits_jv'_j \otimes v''_j$.

Recall \cite{CP} that in the case of $*$-algebras the definition of an
$A$-module algebra includes the following compatibility condition for
involutions:
\begin{equation}\label{*H}
(af)^*=(S(a))^*f^*,\qquad a \in A,\;f \in F.
\end{equation}

Let $U_q \mathfrak{su}_{n,m}$ \index{$U_q \mathfrak{su}_{n,m}$} stands for
the Hopf $*$-algebra $(U_q \mathfrak{sl}_N,*)$ given by
$$
(K_j^{\pm 1})^*=K_j^{\pm 1},\qquad E_j^*=
\begin{cases}
K_jF_j,& j \ne n
\\ -K_jF_j,& j=n
\end{cases},\qquad F_j^*=
\begin{cases}
E_jK_j^{-1},& j \ne n
\\ -E_jK_j^{-1},& j=n
\end{cases},
$$
with $j=1,\ldots,N-1$.
Recall the standard method which was used in \cite{SV2} to equip the space
$$
\mathrm{Pol}(\mathcal{M}at_{m,n})_q
\stackrel{\mathrm{def}}{=}\mathbb{C}[\mathcal{M}at_{m,n}]_q \otimes
\mathbb{C}[\overline{\mathcal{M}at}_{m,n}]_q,
$$
with a structure of $U_q \mathfrak{su}_{n,m}$-module $*$-algebra. The
involution in question allows one, in particular, to introduce the standard
generators $(z_a^\alpha)^*$ of the subalgebra
$\mathbb{C}[\overline{\mathcal{M}at}_{m,n}]_q$. Define the product of
$\varphi_+ \otimes \varphi_-,\,\psi_+ \otimes
\psi_-\:\in\:\mathbb{C}[\mathcal{M}at_{m,n}]_q \otimes
\mathbb{C}[\overline{\mathcal{M}at}_{m,n}]_q$ as follows:
$$
(\varphi_+ \otimes \varphi_-)(\psi_+ \otimes \psi_-)\,=\,m_+ \otimes
m_-\left(\varphi_+ \otimes \check{R}(\varphi_- \otimes \psi_+)\otimes
\psi_-\right).
$$
Here $m_+,\,m_-$ are the multiplications in
$\mathbb{C}[\mathcal{M}at_{m,n}]_q$,
$\mathbb{C}[\overline{\mathcal{M}at}_{m,n}]_q$ respectively, and
\begin{gather}
\check{R}:\mathbb{C}[\overline{\mathcal{M}at}_{m,n}]_q \otimes
\mathbb{C}[\mathcal{M}at_{m,n}]_q \to \mathbb{C}[\mathcal{M}at_{m,n}]_q
\otimes \mathbb{C}[\overline{\mathcal{M}at}_{m,n}]_q,\nonumber
\\ \check{R}=\sigma \cdot R_{\mathbb{C}[\overline{\mathcal{M}at}_{m,n}]_q,
\mathbb{C}[\mathcal{M}at_{m,n}]_q},\label{rhat}
\end{gather}
with $\sigma:a \otimes b \mapsto b \otimes a$. The associativity of the
multiplication in $\mathrm{Pol}(\mathcal{M}at_{m,n})_q$ can be easily
derived from (\ref{braid1}), (\ref{braid2}) by a standard argument
\cite{JS}.
Note that $m_+$, $m_-$, $\check{R}$ are morphisms of $U_q
\mathfrak{sl}_N$-modules. So, $\mathrm{Pol}(\mathcal{M}at_{m,n})_q$ is a
$U_q \mathfrak{sl}_N$-module algebra. We intend to equip
$\mathrm{Pol}(\mathcal{M}at_{m,n})_q$ with an involution. Consider the
antilinear operators $*:V^l \to V^h;\quad *:V^h \to V^l$, which are
determined by the following properties. Firstly, $(v^h)^*=v^l$,
$(v^l)^*=v^h$, and, secondly,
\begin{equation}\label{inv}
(\xi v)^*\,=\,\left(S^{-1}(\xi)\right)^*v^*
\end{equation}
for all $v \in V^h$, (resp. $V^l$), $\xi \in U_q \mathfrak{su}_{n,m}$. This
is certainly equivalent to
$$
\left(\xi v^h \right)^*=\left(S^{-1}(\xi)\right)^*(v^h)^*;\qquad \left(\xi
v^l \right)^*=\left(S^{-1}(\xi)\right)^*(v^l)^*,\qquad \xi \in U_q
\mathfrak{su}_{n,m}.
$$
It follows from the definitions of $V^h$, $V^l$ that the operators as above
are well defined. In particular, (\ref{inv}) can be easily deduced. It also
follows from the relation
$\left(S^{-1}\left(\left(S^{-1}(\xi)\right)^*\right)\right)^*=\xi$ that the
operators are mutually inverse. The duality argument allows one to form the
mutually inverse antihomomorphisms $*:\mathbb{C}[\mathcal{M}at_{m,n}]_q \to
\mathbb{C}[\overline{\mathcal{M}at}_{m,n}]_q ,\quad
*:\mathbb{C}[\overline{\mathcal{M}at}_{m,n}]_q \to
\mathbb{C}[\mathcal{M}at_{m,n}]_q:$
\begin{equation}\label{con1}
f^*(v)\stackrel{\mathrm{def}}{=}\overline{f(v^*)},\qquad v \in
V^l\;\left(\text{resp.}\;V^h \right),\quad f \in \left(V^h
\right)^*\;\left(\text{resp.}\;\left(V^l \right)^*\right).
\end{equation}
Now we are in a position to define the antilinear operator $*$ in
$\mathrm{Pol}(\mathcal{M}at_{m,n})_q$ by
$$(f_+ \otimes f_-)^*\stackrel{\mathrm{def}}{=}f_-^*\otimes f_+^*,$$
for $f_+\in \mathbb{C}[\mathcal{M}at_{m,n}]_q$, $f_-\in
\mathbb{C}[\overline{\mathcal{M}at}_{m,n}]_q$. What remains is to verify
that $*$ equips $\mathrm{Pol}(\mathcal{M}at_{m,n})_q$ with a structure of
$U_q \mathfrak{su}_{n,m}$-module algebra. For that, the reader is referred
to \cite[section 8]{SV2}.

We identify $\mathbb{C}[\mathcal{M}at_{m,n}]_q$ with its image under the
embedding $\mathbb{C}[\mathcal{M}at_{m,n}]_q \hookrightarrow
\mathrm{Pol}(\mathcal{M}at_{m,n})_q$, $f \mapsto f \otimes 1$, and
$\mathbb{C}[\overline{\mathcal{M}at}_{m,n}]_q$ with its image under the
embedding $\mathbb{C}[\overline{\mathcal{M}at}_{m,n}]_q \hookrightarrow
\mathrm{Pol}(\mathcal{M}at_{m,n})_q$, $f \mapsto 1 \otimes f$.
It follows from proposition \ref{gen} that $\{z_a^\alpha \}$,
$\alpha=1,\ldots,m$, $a=1,\ldots,n$, generate the $*$-algebra
$\mathrm{Pol}(\mathcal{M}at_{m,n})_q$, and the complete list of relations
consists of (\ref{zaa1}) -- (\ref{zaa3}), together with the following one:

\begin{equation}\label{zsz}
\left(z_b^\beta \right)^*z_a^\alpha=m \check{R}\left(\left(z_b^\beta
\right)^*\otimes z_a^\alpha \right)=m \sigma
R_{\mathbb{C}[\overline{\mathcal{M}at}_{m,n}]_q
\,\mathbb{C}[\mathcal{M}at_{m,n}]_q}\left(\left(z_b^\beta \right)^*\otimes
z_a^\alpha \right),
\end{equation}
with $m:\mathrm{Pol}(\mathcal{M}at_{m,n})_q^{\otimes 2}\to
\mathrm{Pol}(\mathcal{M}at_{m,n})_q$, $m:f_1 \otimes f_2 \mapsto f_1f_2$
being the multiplication in $\mathrm{Pol}(\mathcal{M}at_{m,n})_q$.

Simplify the expression $R_{\mathbb{C}[\overline{\mathcal{M}at}_{m,n}]_q
\,\mathbb{C}[\mathcal{M}at_{m,n}]_q}\left(\left(z_b^\beta \right)^*\otimes
z_a^\alpha \right)$ in (\ref{zsz}). Denote by $U_q \mathfrak{su}_n \otimes
U_q \mathfrak{su}_m$ \index{$U_q \mathfrak{su}_n \otimes U_q
\mathfrak{su}_m$} the subalgebra of the Hopf $*$-algebra $U_q
\mathfrak{su}_{n,m}$ generated by $E_j$, $F_j$, $K_j$, $K_j^{-1}$ with $j
\ne n$. Now an application of proposition \ref{zaalph} makes it easy to
prove the following

\medskip

\begin{lemma}\label{spi}
The sesquilinear form in $\mathbb{C}[\mathcal{M}at_{m,n}]_{q,1}$ given by
$\left(z_a^\alpha,z_b^\beta \right)=\delta_{ab}\delta^{\alpha \beta}$,
$a,b=1,\ldots,n$, $\alpha,\beta=1,\ldots,m$, is $U_q \mathfrak{su}_n \otimes
U_q \mathfrak{su}_m$-invariant:
$$
\left(\xi z_a^\alpha,z_b^\beta \right)=\left(z_a^\alpha,\xi^*z_b^\beta
\right),\qquad \xi \in U_q \mathfrak{su}_n \otimes U_q \mathfrak{su}_m,\quad
a,b=1,\ldots,n,\quad \alpha,\beta=1,\ldots,m.
$$
\end{lemma}

\medskip

\begin{corollary}\label{mu}
The linear functional $\mu$ on
$\mathbb{C}[\overline{\mathcal{M}at}_{m,n}]_{q,-1}\otimes
\mathbb{C}[\mathcal{M}at_{m,n}]_{q,1}$ given by $\mu \left(\left(z_b^\beta
\right)^*\otimes z_a^\alpha \right)=\delta_{ab} \delta^{\alpha \beta}$, is
invariant:
\begin{multline*}
\mu \left(\xi \left(\left(z_b^\beta \right)^*\otimes z_a^\alpha
\right)\right)=\varepsilon(\xi)\mu \left(\left(z_b^\beta \right)^*\otimes
z_a^\alpha \right),
\\ \xi \in U_q \mathfrak{su}_n \otimes U_q \mathfrak{su}_m,\quad
a,b=1,\ldots,n,\quad \alpha,\beta=1,\ldots,m.
\end{multline*}
\end{corollary}

\smallskip

{\bf Proof.} Let $L=\mathbb{C}[\mathcal{M}at_{m,n}]_{q,1}$. Consider the
{\it antimodule} $\overline{L}$ which is still $L$ as an Abelian group, but
the actions of the ground field and $U_q \mathfrak{su}_n \otimes U_q
\mathfrak{su}_m$ are given by $(\lambda,v)\mapsto \overline{\lambda}v$,
$(\xi,v)\mapsto S(\xi)^*v$, $\xi \in U_q \mathfrak{su}_n \otimes U_q
\mathfrak{su}_m$, $v \in L$. It follows from lemma \ref{spi} that the linear
functional $\overline{L}\otimes L \to \mathbb{C}(q^{1/s})$, corresponding to
the sesquilinear form in $L$, is invariant. \hfill $\square$

\medskip

Let $L'=\mathbb{C}[\overline{\mathcal{M}at}_{m,n}]_{q,-1}$,
$L''=\mathbb{C}[\mathcal{M}at_{m,n}]_{q,1}$, and $R_{L'L''}$ is the linear
operator in $L'\otimes L''$ given by the action of the universal R-matrix of
the Hopf algebra $U_q \mathfrak{sl}_n \otimes U_q \mathfrak{sl}_m \subset
U_q \mathfrak{sl}_N$, determined by a formula similar to (\ref{rmmf}).

\medskip

\begin{lemma}\label{c1c2}
In $\mathrm{Pol}(\mathcal{M}at_{m,n})_q$, for all $a,b=1,\ldots,n$,
$\alpha,\beta=1,\ldots,m$,
\begin{multline}\label{c1c2_}
R_{\mathbb{C}[\overline{\mathcal{M}at}_{m,n}]_q
\,\mathbb{C}[\mathcal{M}at_{m,n}]_q}\left(\left(z_b^\beta \right)^*\otimes
z_a^\alpha \right)=
\\ ={\tt const}_1 \cdot R_{L'L''}\left(\left(z_b^\beta
\right)^*\otimes z_a^\alpha \right)+{\tt const}_2 \cdot
\delta_{ab}\delta^{\alpha \beta},
\end{multline}
with ${\tt const}_1$ and ${\tt const}_2$ being independent of $a$, $b$,
$\alpha$, $\beta$.
\end{lemma}

\smallskip

{\bf Proof.} Reduce the left hand side of (\ref{c1c2_}) modulo
$\mathbb{C}[\overline{\mathcal{M}at}_{m,n}]_{q,0}\otimes
\mathbb{C}[\mathcal{M}at_{m,n}]_{q,0}$, using (\ref{rmmf}). The 'redundant'
exponential multiples in the left hand side of the resulting identity can be
omitted since
$$
\exp_{q^2}\left((q^{-1}-q)E_{\beta_j}\otimes
F_{\beta_j}\right)\left(\left(z_b^\beta \right)^*\otimes
z_a^\alpha\right)=\left(z_b^\beta \right)^*\otimes z_a^\alpha
$$
for all $\alpha,\beta=1,\ldots,m$, $a,b=1,\ldots,n$,
$j>\dfrac{m(m-1)}2+\dfrac{n(n-1)}2$. What remains is to compare the multiple
$q^{-t_0}$ related to the Hopf algebra $U_q \mathfrak{sl}_N$ to a similar
multiple related to the Hopf subalgebra $U_q \mathfrak{sl}_n \otimes U_q
\mathfrak{sl}_m$. It follows from lemma \ref{orth} and the description of
$t_0$ in terms of the orthogonal basis of the Cartan subalgebra
(\ref{3descr}) that their actions on the subspace
$\mathbb{C}[\overline{\mathcal{M}at}_{m,n}]_{q,-1}\otimes
\mathbb{C}[\mathcal{M}at_{m,n}]_{q,1}$ differ only by a constant multiple.
This implies the existence of such element ${\tt const}_1$ of the ground
field that for all $a$, $b$, $\alpha$, $\beta$ one has
\begin{multline*}
R_{\mathbb{C}[\overline{\mathcal{M}at}_{m,n}]_q
\,\mathbb{C}[\mathcal{M}at_{m,n}]_q}\left(\left(z_b^\beta \right)^*\otimes
z_a^\alpha \right)-{\tt const}_1 \cdot R_{L'L''}\left(\left(z_b^\beta
\right)^*\otimes z_a^\alpha \right)\in
\\ \in \mathbb{C}[\overline{\mathcal{M}at}_{m,n}]_{q,0}\otimes
\mathbb{C}[\mathcal{M}at_{m,n}]_{q,0}.
\end{multline*}
Thus we get a linear functional $l$ on
$\mathbb{C}[\overline{\mathcal{M}at}_{m,n}]_{q,-1}\otimes
\mathbb{C}[\mathcal{M}at_{m,n}]_{q,1}$ since
$$
\dim(\mathbb{C}[\overline{\mathcal{M}at}_{m,n}]_{q,0})=
\dim(\mathbb{C}[\mathcal{M}at_{m,n}]_{q,0})=1.
$$
Clearly, the linear maps $\sigma \cdot
R_{\mathbb{C}[\overline{\mathcal{M}at}_{m,n}]_q
\,\mathbb{C}[\mathcal{M}at_{m,n}]_q}$ and $\sigma \cdot R_{L'L''}$ are
morphisms of $U_q \mathfrak{sl}_n \otimes U_q \mathfrak{sl}_m$-modules. So
the linear functional $l$ is $U_q \mathfrak{sl}_n \otimes U_q
\mathfrak{sl}_m$-invariant. What remains is to apply corollary \ref{mu},
together with the fact that the subspace of $U_q \mathfrak{sl}_n \otimes U_q
\mathfrak{sl}_m$-invariant functionals
$\mathbb{C}[\overline{\mathcal{M}at}_{m,n}]_{q,-1}\otimes
\mathbb{C}[\mathcal{M}at_{m,n}]_{q,1}\to \mathbb{C}$ is
one-dimensional.\footnote{In fact, the dimensions of isotypic components of
the $U_q \mathfrak{sl}_n \otimes U_q \mathfrak{sl}_m$-module
$\mathbb{C}[\overline{\mathcal{M}at}_{m,n}]_{q,-1}\otimes
\mathbb{C}[\mathcal{M}at_{m,n}]_{q,1}$ are the same just as in the case
$q=1$.} \hfill $\square$

\medskip

We need an explicit form of the operator $R_{L'L''}$. Let $*:U \to
\overline{U}$, $*:V \to \overline{V}$, be the identical maps from the above
$U_q \mathfrak{sl}_n$-module $U$ and $U_q \mathfrak{sl}_m$-module $V$ onto
the associated antimodules. Let $R_{\overline{U}U}$, $R_{\overline{V}V}$
stand for the operators in $\overline{U}\otimes U$, $\overline{V}\otimes V$
respectively, given by the actions of the universal R-matrices of the Hopf
algebras $U_q \mathfrak{sl}_n$ and $U_q \mathfrak{sl}_m$. The following
result is well known; we reproduce its proof here for the reader's
convenience.

\medskip

\begin{lemma}\label{c'c''}
For all $a,b=1,\ldots,n$, $\alpha,\beta=1,\ldots,m$,
\begin{eqnarray*}
R_{\overline{U}U}(u_b^*\otimes u_a)&=&{\tt const}'\cdot
\begin{cases}q^{-1}u_b^*\otimes u_a, & a \ne b
\\ u_a^*\otimes u_a-(q^{-2}-1)\sum \limits_{k>a}u_k^* \otimes u_k,&a=b
\end{cases},
\\ R_{\overline{V}V}((v^\beta)^*\otimes v^\alpha)&=&{\tt const}''\cdot
\begin{cases}q^{-1}\left(v^\beta \right)^*\otimes v^\alpha,& \alpha \ne \beta
\\ \left(v^\alpha \right)^*\otimes v^\alpha-(q^{-2}-1)\sum
\limits_{k>\alpha}\left(v^k \right)^*\otimes v^k,& \alpha=\beta
\end{cases}.
\end{eqnarray*}
with ${\tt const}'$, ${\tt const}''$ being independent of $a$, $b$,
$\alpha$, $\beta$.
\end{lemma}

\smallskip

{\bf Proof.} It suffices to prove the first identity. Consider the linear
operator $\sigma \cdot R_{\overline{U}U}:\overline{U}\otimes U \to{U}\otimes
\overline{U}$, with $\sigma$ being the flip of tensor multiples. This
operator is a morphism of $U_q \mathfrak{sl}_n$-modules. Besides, it follows
from (\ref{rmmf}) that $\sigma \cdot R_{\overline{U}U}(u_n^* \otimes
u_n)={\tt const}'\cdot u_n \otimes u_n^*$ since $u_n$ is the lowest weight
vector of the $U_q \mathfrak{sl}_n$-module $U$. On the other hand, it is
well known (see, for example, \cite{SoV1}) that the composition of $\sigma$
with the operator defined by the right hand side of the first identity in
the statement of our lemma possesses the same properties. What remains is to
use the fact that each morphism of $U_q \mathfrak{sl}_n$-modules
$\overline{U}\otimes U \to{U}\otimes \overline{U}$ which annihilates $u_n^*
\otimes u_n$, is identically zero (this vector does not belong to any of the
two simple components of the $U_q \mathfrak{sl}_n$-module
$\overline{U}\otimes U$, and hence it generates this module). \hfill
$\square$

\medskip

Lemmas \ref{c1c2}, \ref{c'c''} allow one to deduce all the relations between
$(z_b^\beta)^*$, $z_a^\alpha$ up to two constants. These will be computed by
means of the following

\medskip

\begin{lemma}\label{znm}
$R_{\mathbb{C}[\overline{\mathcal{M}at}_{m,n}]_q
\,\mathbb{C}[\mathcal{M}at_{m,n}]_q}((z_n^m)^*\otimes
z_n^m)=q^2(z_n^m)^*\otimes z_n^m+1-q^2$.
\end{lemma}

\smallskip

{\bf Proof.} We are about to apply the explicit formula (\ref{rmmf}) for the
universal R-matrix. Prove that
$$H_jz_n^m=
\begin{cases}
2z_n^m,& j=n
\\ -z_n^m,& |j-n|=1
\\ 0,& \mathrm{otherwise}
\end{cases}.
$$
The two latter relations follow from the definitions of $z_a^\alpha$. The
first relation follows from $H_0z_n^m=2z_n^m$:
$$2z_n^m=\frac2{m+n}(-m(n-1)-n(m-1))z_n^m+\frac{2mn}{m+n}H_nz_n^m.$$
Hence $z_n^m$, $(z_n^m)^*$ are weight vectors whose weights are $\alpha_n$,
$-\alpha_n$ respectively. Thus, we have
$$
t_0((z_n^m)^*\otimes z_n^m)=(-\alpha_n,\alpha_n)(z_n^m)^*\otimes
z_n^m=-2(z_n^m)^*\otimes z_n^m.
$$
Observe that the only $q$-exponent, which survives in (\ref{rmmf}), is
$\exp_{q^2}\left((q^{-1}-q)E_{\beta_M}\otimes F_{\beta_M}\right)$. Of
course, $\beta_M=\alpha_n$, and it is not difficult to prove that
$E_{\beta_M}\otimes F_{\beta_M}=\mathrm{const}\cdot E_n\otimes F_n$. The
constant multiplier in the latter equality is equal to $1$ since
$$\left(F_{\beta_M}, E_{\beta_M}\right)=
\left(F_n, E_n\right)=\frac1{q^{-1}-q}$$ with respect to the well known
pairing $(\,,\,):U_q\mathfrak{b}^-\times U_q\mathfrak{b}^+\to\mathbb{C}$
\cite[Chapter 6]{Jant}.\footnote{This pairing is used in constructing of the
canonical homomorphism $D(U_q\mathfrak{b}^+)\to U_q\mathfrak{sl}_N$ where
$D(U_q\mathfrak{b}^+)$ is the quantum double of the Hopf algebra
$U_q\mathfrak{b}^+$ \cite{Dr1}.} Hence
$$\exp_{q^2}\left((q^{-1}-q)E_{\beta_M}\otimes
F_{\beta_M}\right)=\exp_{q^2}\left((q^{-1}-q)E_n\otimes F_n\right),$$ and
what remains is to use the detailed calculations for the case $m=n=1$ given
in \cite{SV2}.

\hfill $\square$

\medskip

\begin{corollary}\label{znmcr}
$(z_n^m)^*z_n^m=q^2z_n^m(z_n^m)^*+1-q^2$.
\end{corollary}

\medskip

Recall that it was our convention at the beginning of this section that our
ground field is $\mathbb{C}(q^{1/s})$. Now let us transfer literally the
definitions of $\mathbb{C}[\mathcal{M}at_{m,n}]_q$ and
$\mathrm{Pol}(\mathcal{M}at_{m,n})_q$ onto the case of $\mathbb{C}$ as a
ground field. In what follows, this will be our context. The result of
section \ref{PBW} implies
\begin{equation}\label{dimm}
 \dim
\mathbb{C}[\mathcal{M}at_{m,n}]_{q,i}=\dbinom{mn+i-1}{i}.
\end{equation}
As a consequence of lemmas \ref{c1c2} -- \ref{znm} and corollary \ref{hom}
one also has

\medskip

\begin{proposition}\label{hom1}
There exists a unique homomorphism of $*$-algebras
\begin{equation}\label{isomor}
\mathrm{Pol}(\mathrm{Mat}_{m,n})_q \to \mathrm{Pol}(\mathcal{M}at_{m,n})_q
\end{equation}
such that $z_a^\alpha \mapsto z_a^\alpha$, $a=1,\ldots,n$,
$\alpha=1,\ldots,m$.
\end{proposition}

\medskip

{\sc Remark.} It will be proved in section \ref{pol} that this homomorphism
appears to be an isomorphism.

\bigskip

\section{The algebras \boldmath $\mathbb{C}[SL_N]_{q,t}$ and
$\mathbb{C}[\widetilde{G}]_{q,x}$}\label{local}

Recall that the $U_q \mathfrak{sl}_N$-module algebra $\mathbb{C}[SL_N]_q$ is
a domain, and its element $t=t_{\{1,2,\dots,m \}\{n+1,n+2,\dots,N
\}}^{\wedge m}$ quasi-commutes with all the generators $t_{ij}$. Let
$\mathbb{C}[SL_N]_{q,t}$ stand for the localization of the algebra
$\mathbb{C}[SL_N]_q$ with respect to the multiplicative set generated by
$t$.

\begin{proposition}\label{texext}
There exists a unique extension of the structure of $U_q
\mathfrak{sl}_N$-module algebra from $\mathbb{C}[SL_N]_q$ onto
$\mathbb{C}[SL_N]_{q,t}$.
\end{proposition}

\smallskip

{\bf Proof.} The uniqueness of the extension is obvious. We are going to
construct such extension by applying the following statement.

\medskip

\begin{lemma}
For every $\xi \in U_q \mathfrak{sl}_N$, $f \in \mathbb{C}[SL_N]_q$, there
exists a unique Laurent polynomial $p_{\xi,f}(\lambda)$ with coefficients
from $\mathbb{C}[SL_N]_{q,t}$ such that
$$
p_{\xi,f}\left(q^l \right)=\xi \left(f \cdot t^l \right)\cdot t^{-l},\qquad
l \in \mathbb{Z}_+.
$$
\end{lemma}

\smallskip

{\bf Proof.} Our statement follows from the definition of a $U_q
\mathfrak{sl}_N$-module algebra structure in $\mathbb{C}[SL_N]_q$
(\ref{EFt}), (\ref{Kt}), the definition of a comultiplication in $U_q
\mathfrak{sl}_N$ (\ref{comult}), and (\ref{ttij}). \hfill $\square$

\medskip

Turn back to the proof of proposition \ref{texext}. We can use the same
Laurent polynomials for defining $\xi \left(f \cdot t^l \right)$ for $f \in
\mathbb{C}[SL_N]_q$ and all integers $l$:
$$
\xi \left(f \cdot t^l \right)\stackrel{\mathrm{def}}{=}p_{\xi,f}\left(q^l
\right)t^l.
$$
Of course, we need firstly to verify that the map $U_q \mathfrak{sl}_N
\times \mathbb{C}[SL_N]_{q,t}\to \mathbb{C}[SL_N]_{q,t}$, $\xi \times f
\mapsto \xi(f)$ as above is well defined, and secondly that we obtain this
way a structure of $U_q \mathfrak{sl}_N$-module algebra. The first item is
equivalent to
$$
p_{\xi,f}\left(q^{a+l}\right)t^a=p_{\xi,f \cdot t^a}\left(q^l \right),\qquad
\xi \in U_q \mathfrak{sl}_N,\;f \in \mathbb{C}[SL_N]_q,\;a \in
\mathbb{Z}_+,\;l \in \mathbb{Z}.
$$
This relation is obvious for $l \in \mathbb{Z}_+$, hence it is valid for all
integers $l$ due to the well known uniqueness theorem for the Laurent
polynomials.

For the second item, we have to prove some identities for $\xi \in U_q
\mathfrak{sl}_N$, $f_1 \cdot t^l$, $f_2 \cdot t^l$, $f_1,f_2 \in
\mathbb{C}[SL_N]_q$, $l \in \mathbb{Z}$. Observe that the left and right
hand sides of those identities (up to multiplying by the same powers of $t$)
are just Laurent polynomials of the indeterminate $\lambda=q^l$. So, it
suffices to prove them for $l \in \mathbb{Z}_+$ due to the same uniqueness
theorem for Laurent polynomials. On the other hand, at all $l \in
\mathbb{Z}_+$ one can deduce these identities from the fact that
$\mathbb{C}[SL_N]_q$ is a $U_q \mathfrak{sl}_N$-module algebra. \hfill
$\square$

\medskip

A more general but less elementary approach to proving statements like
proposition \ref{texext} have been obtained in a recent work by V. Lunts and
A. Rosenberg \cite{LR}.

The following result is due to M. Noumi \cite{N} for the case $m=n=2$. It
will be also refined in the sequel. Recall the notation $J_{a
\alpha}=\{n+1,n+2,\ldots,N \}\setminus \{N+1-\alpha \}\cup \{a \}$
\index{$J_{a \alpha}$} as in the statement of proposition \ref{emb}.

\medskip

\begin{proposition}\label{emb1}
The map $i:z_a^\alpha \mapsto t^{-1}\cdot t_{\{1,2,\ldots,m \}J_{a
\alpha}}^{\wedge m}$, $\alpha=1,\ldots,m$, $a=1,\ldots,n$, admits a unique
extension up to a homomorphism of algebras
$i:\mathbb{C}[\mathrm{Mat}_{m,n}]_q \to \mathbb{C}[SL_N]_{q,t}$.
\end{proposition}

\smallskip

{\sc Remark.} It will be shown later that in fact $i$ is an embedding (see
section \ref{ii}).

\smallskip

{\bf Proof} of proposition \ref{emb1}. The uniqueness of the extension is
obvious.

Let $\mathbb{C}[\mathrm{Mat}_{m,N}]_q$ be the algebra defined by its
generators $\{t_{\alpha a}\}$, $\alpha=1,2,\ldots,m$, $a=1,2,\ldots,N$, and
the relations (\ref{taa1}) -- (\ref{taa3}), together with
$t=t_{\{1,2,\dots,m \}\{n+1,n+2,\dots,N \}}^{\wedge m}$. Consider the
localization $\mathbb{C}[\mathrm{Mat}_{m,N}]_{q,t}$ of
$\mathbb{C}[\mathrm{Mat}_{m,N}]_q$ with respect to the multiplicatively
closed set $t^\mathbb{N}$. It suffices to prove that the map
$$
i:z_a^\alpha \mapsto t^{-1}\cdot t_{\{1,2,\ldots,m \}J_{a \alpha}}^{\wedge
m},\qquad \alpha=1,\ldots,m,\quad a=1,\ldots,n
$$
admits an extension up to a homomorphism of algebras
$i:\mathbb{C}[\mathrm{Mat}_{m,n}]_q \to
\mathbb{C}[\mathrm{Mat}_{m,N}]_{q,t}$. Consider an embedding of $U_q
\mathfrak{sl}_N$-module algebras $i':\mathbb{C}[\mathrm{Mat}_{m,N}]_q
\hookrightarrow (U_q \mathfrak{sl}_N)^*$, which is a composition of the
embedding $\mathbb{C}[\mathrm{Mat}_{m,N}]_q \hookrightarrow
\mathbb{C}[SL_N]_q$, $t_{\alpha,a}\mapsto t_{\alpha+n,a}$,
$\alpha=1,2,\ldots,m$, $a=1,2,\ldots,N$, and the canonical embedding
$\mathbb{C}[SL_N]_q \hookrightarrow (U_q \mathfrak{sl}_N)^*$. One can use
the same argument as in the proof of proposition \ref{texext} to extend the
structure of $U_q \mathfrak{sl}_N$-module algebra and the embedding $i'$
onto the localization $\mathbb{C}[\mathrm{Mat}_{m,N}]_{q,t}$ of
$\mathbb{C}[\mathrm{Mat}_{m,N}]_q$.

Consider the embedding $i'':\mathbb{C}[\mathcal{M}at_{m,n}]_q
\hookrightarrow (U_q \mathfrak{sl}_N)^*$, derived by a duality from the onto
morphism of coalgebras
$$
j:U_q \mathfrak{sl}_N \to V^h, \qquad j:\xi \mapsto S(\xi)v^h, \quad \xi \in
U_q \mathfrak{sl}_N,
$$
with $v^h$ being the generator of the $U_q \mathfrak{sl}_N$-module $V^h$.
Use the embedding $i''$ to get a composition of $i''$ and the homomorphism
$\mathbb{C}[\mathrm{Mat}_{m,n}]_q \to \mathbb{C}[\mathcal{M}at_{m,n}]_q$
(see section \ref{da}) to obtain a homomorphism of algebras:
\begin{equation}\label{k}
i''':\mathbb{C}[\mathrm{Mat}_{m,n}]_q \to(U_q \mathfrak{sl}_N)^*.
\end{equation}
One can easily apply proposition \ref{zaalph} to prove that the statement of
proposition \ref{emb1} reduces to

\medskip

\begin{lemma}\label{i'i''}
$$
i'\left(t^{-1}t_{\{1,2,\ldots,m \}J}^{\wedge m}\right)\subset
i''(\mathbb{C}[\mathcal{M}at_{m,n}]_q),\qquad \mathrm{card}(J)=m.
$$
\end{lemma}

\smallskip

{\bf Proof.} It suffices to establish that $i'\left(t^{-1}t_{\{1,2,\ldots,m
\}J}^{\wedge m}\right)$ are orthogonal to the kernel of $j$ with respect to
the above pairing. Let us agree not to distinguish between the $U_q
\mathfrak{sl}_N$-module algebras $\mathbb{C}[\mathcal{M}at_{m,n}]_q$,
$\mathbb{C}[\mathrm{Mat}_{m,N}]_{q,t}$, $\mathbb{C}[SL_N]_q$ and their
images in $(U_q \mathfrak{sl}_N)^*$. What remains now is to prove that for
$\mathrm{card}(J)=m$, $\xi \in U_q \mathfrak{sl}_N$,
\begin{align*}
& \left \langle t^{-1}t_{\{1,2,\ldots,m \}J}^{\wedge m},(K_i^{\pm 1}-1)\xi
\right \rangle=0,&& i=1,2,\ldots,N-1,
\\ & \left \langle t^{-1}t_{\{1,2,\ldots,m
\}J}^{\wedge m},E_i \xi \right \rangle=0,&& i=1,2,\ldots,N-1,
\\ & \left \langle t^{-1}t_{\{1,2,\ldots,m \}J}^{\wedge m},F_j \xi
\right \rangle=0,&& j=1,2,\ldots,n-1,n+1,\ldots,N-1.
\end{align*}
These follow from the more general relations
\begin{alignat*}{3}
& \left \langle t^k t_{\{1,2,\ldots,m \}J}^{\wedge m},(K_i^{\pm 1}-1)\xi
\right \rangle=&&\left(q^{\mp(k+1)}-1 \right)\delta_{im}\left \langle t^k
t_{\{1,2,\ldots,m \}J}^{\wedge m},\xi \right \rangle,
\\ & \left \langle t^k t_{\{1,2,\ldots,m \}J}^{\wedge m},E_i \xi
\right \rangle=0,&& i=1,2,\ldots,N-1,
\\ & \left \langle t^k t_{\{1,2,\ldots,m \}J}^{\wedge m},F_j \xi \right
\rangle=0,&& j=1,2,\ldots,n-1,n+1,\ldots,N-1,
\end{alignat*}
for $\xi \in U_q \mathfrak{sl}_N$. In proving these latter relations, one
can restrict matters to the case $k \in \mathbb{Z}_+$ by using the
techniques related to Laurent polynomials. Let
$\widetilde{t}=t_{\{n+1,n+2,\dots,N \}\{n+1,n+2,\dots,N \}}^{\wedge m}$.
What remains to prove now is that for all $\xi \in U_q \mathfrak{sl}_N$, $k
\in \mathbb{Z}_+$
\begin{alignat*}{3}
& \left \langle \widetilde{t}^{\,k}t_{\{n+1,n+2,\ldots,N \}J}^{\wedge
m},(K_i^{\pm 1}-1)\xi \right \rangle=&&\left(q^{\mp(k+1)}-1 \right
)\delta_{im}\left \langle \widetilde{t}^{\,k}t_{\{1,2,\ldots,m \}J}^{\wedge
m},\xi \right \rangle,
\\ & \left \langle \widetilde{t}^{\,k}t_{\{n+1,n+2,\ldots,N
\}J}^{\wedge m},E_i \xi \right \rangle=0,&& i=1,2,\ldots,N-1,
\\ & \left \langle \widetilde{t}^{\,k}t_{\{n+1,n+2,\ldots,N
\}J}^{\wedge
m},F_j \eta \right \rangle=0,&& j=1,2,\ldots,n-1,n+1,\ldots,N-1. \qquad
\square
\end{alignat*}

\medskip

{\sc Remark.} The proof of proposition \ref{emb1} involves a construction of
the embedding $i'^{-1}i''$ of $U_q \mathfrak{sl}_N$-module algebras
$\mathbb{C}[\mathcal{M}at_{m,n}]_q\to\mathbb{C}[\mathrm{Mat}_{m,N}]_{q,t}$.
Hence, the map
$$
z_a^\alpha \mapsto t^{-1}\cdot t_{\{1,2,\ldots,m \}J_{a \alpha}}^{\wedge
m},\qquad \alpha=1,\ldots,m,\quad a=1,\ldots,n.
$$
admits an extension up to an embedding of $U_q \mathfrak{sl}_N$-module
algebras $\mathcal{I}:\mathbb{C}[\mathcal{M}at_{m,n}]_q\to
\mathbb{C}[SL_N]_{q,t}$.


\medskip

\begin{lemma}\label{Noumi}
For all $1 \le a<b \le n$, $1 \le \alpha<\beta \le m$,
$$
\mathcal{I}:z_a^\alpha z_b^\beta-qz_a^\beta z_b^\alpha \mapsto t^{-1}\cdot
t_{\{1,2,\ldots,m
\}\{a,b,\ldots,\widehat{N+1-\beta},\ldots,\widehat{N+1-\alpha},\ldots,N
\}}^{\wedge m}.
$$
\end{lemma}

\smallskip

{\bf Proof.} In the same way as in the proof of proposition \ref{emb1}, one
can establish that $t^{-1}\cdot t_{\{1,2,\ldots,m
\}\{a,b,\ldots,\widehat{N+1-\beta},\ldots,\widehat{N+1-\alpha},\ldots,N
\}}^{\wedge m}\in \mathcal{I}(\mathbb{C}[\mathcal{M}at_{m,n}]_q)$. This is a
weight vector of the $U_q \mathfrak{sl}_N$-module
$\mathcal{I}(\mathbb{C}[\mathcal{M}at_{m,n}]_q)$. A computation of the
weight yields
$$
\mathcal{I}(c_1z_a^\beta z_b^\alpha+c_2z_a^\alpha z_b^\beta)=t^{-1}\cdot
t_{\{1,2,\ldots,m
\}\{a,b,\ldots,\widehat{N+1-\beta},\ldots,\widehat{N+1-\alpha},\ldots,N
\}}^{\wedge m},
$$
with $c_1$, $c_2\in \mathbb{C}$. When computing the constants $c_1$, $c_2$,
one can restrict oneself to the special case \footnote{The general case is
derivable by observing that $\mathcal{I}$ is a morphism of $U_q
\mathfrak{sl}_N$-modules.} $a=n-1$, $b=n$, $\alpha=m-1$, $\beta=m$. Even
more, one can stick to the case $m=n=2$ due to the homomorphism:
$$
t_{i,j}\mapsto
\begin{cases}
t_{i,j-n+2},& i \le 2 \quad \& \quad n-1 \le j \le n+2,
\\ 1,& i>2 \quad \& \quad j=i+n,
\\ 0,& \text{otherwise}.
\end{cases}
$$
In the special case $m=n=2$ the result in question is accessible via a
direct calculation \cite{N}. \hfill $\square$

\medskip

A proof of a more general statement is presented in section \ref{pry}.

The lemma \ref{Noumi}, together with proposition \ref{zaalph}, allow one to
get a description of the $U_q \mathfrak{sl}_N$-module algebra structure on
$\mathbb{C}[\mathcal{M}at_{m,n}]_q$.

\begin{corollary}\label{uqslmat}
For ${a=1,\ldots,n;\alpha=1,\ldots,m}$
\begin{gather*}
K_n^{\pm 1}z_a^\alpha=
\begin{cases}
q^{\pm 2}z_a^\alpha,&a=n \;\&\;\alpha=m
\\ q^{\pm 1}z_a^\alpha,&a=n \;\&\;\alpha \ne m \quad \mathrm{or}\quad a \ne
n \;\&\; \alpha=m
\\ z_a^\alpha,&\mathrm{otherwise}
\end{cases}
\\ F_nz_a^\alpha=q^{1/2}\cdot
\begin{cases}
1,& a=n \;\& \;\alpha=m
\\ 0,&\mathrm{otherwise}
\end{cases}\qquad
E_nz_a^\alpha=-q^{1/2}\cdot
\begin{cases}
q^{-1}z_a^mz_n^\alpha,&a \ne n \;\&\;\alpha \ne m
\\ (z_n^m)^2,& a=n \;\&\;\alpha=m
\\ z_n^mz_a^{\alpha},&\mathrm{otherwise}
\end{cases}
\end{gather*}
and with $k \ne n$
\begin{align*}
K_k^{\pm 1}z_a^\alpha&=
\begin{cases}
q^{\pm 1}z_a^\alpha,& k<n \;\&\;a=k \quad \mathrm{or}\quad k>n
\;\&\;\alpha=N-k
\\ q^{\mp 1}z_a^\alpha,& k<n \;\&\;a=k+1 \quad \mathrm{or}\quad k>n
\;\&\;\alpha=N-k+1
\\ z_a^\alpha,&\mathrm{ otherwise}
\end{cases}
\\ F_kz_a^\alpha&=q^{1/2}\cdot
\begin{cases}
z_{a+1}^\alpha,& k<n \;\&\;a=k
\\ z_a^{\alpha+1},& k>n \;\&\;\alpha=N-k
\\ 0,&\mathrm{otherwise}
\end{cases},
\\ E_kz_a^\alpha&=q^{-1/2}\cdot
\begin{cases}
z_{a-1}^\alpha,& k<n \;\&\; a=k+1
\\ z_a^{\alpha-1},& k>n \;\&\;\alpha=N-k+1
\\ 0,& \mathrm{otherwise}
\end{cases}.
\end{align*}
\end{corollary}

\bigskip

In sections \ref{qphs}, \ref{da} the algebras $\mathbb{C}[SL_N]_q$ and $U_q
\mathfrak{sl}_N$ were equipped with involutions. Thus we got $*$-algebras
$\mathbb{C}[\widetilde{G}]_q$ and $U_q \mathfrak{su}_{n,m}$. Recall that
$\mathbb{C}[SL_N]_q$ is a $U_q \mathfrak{sl}_N$-module algebra. It is easy
to prove that the involutions in question agree in such a way that
$\mathbb{C}[\widetilde{G}]_q=(\mathbb{C}[SL_N]_q,*)$ is a $U_q
\mathfrak{su}_{n,m}$-module algebra\footnote{It suffices to use a similar
result for the quantum group $SU_N$ and the quantum universal enveloping
algebra $U_q \mathfrak{su}_N$, as the compact and the non-compact
involutions on the generators differ only by sign change.}.

In section \ref{qphs} an element $x=tt^*$ and the localization
$\mathbb{C}[\widetilde{G}]_{q,x}$ were introduced. An argument similar to
that used in the proof of proposition \ref{texext} (with the reference to
(\ref{xtij}) instead of (\ref{ttij})) allows to get

\begin{proposition}\label{texext1}
There exists a unique extension of the structure of $U_q
\mathfrak{su}_{n,m}$-module algebra from $\mathbb{C}[\widetilde{G}]_q$ onto
$\mathbb{C}[\widetilde{G}]_{q,x}$.
\end{proposition}

\begin{proposition}\label{emb3}
The map
\begin{equation}\label{cI}
i:z_a^\alpha \mapsto t^{-1}t_{\{1,2,\ldots,m \}J_{a \alpha}}^{\wedge m}
\end{equation}
with $J_{a \alpha}=\{n+1,n+2,\ldots,N \}\setminus \{N+1-\alpha \}\cup \{a
\}$, is uniquely extendable up to a homomorphism of $*$-algebras
$i:\mathrm{Pol}(\mathrm{Mat}_{m,n})_q \to \mathbb{C}[\widetilde{G}]_{q,x}$.
\end{proposition}

{\sc Remark.} It will be shown later that in fact $i$ is an embedding (see
the conclusion remark in section \ref{ii}).

{\bf Proof} of proposition \ref{emb3}. The uniqueness of the extension is
obvious. The existence follows from a construction of a homomorphism of
$U_q\mathfrak{su}_{n,m}$-module algebras
\begin{equation}\label{emb4}
\mathrm{Pol}(\mathcal{M}at_{m,n})_q \to
\mathbb{C}[\widetilde{G}]_{q,x},\qquad z_a^\alpha \mapsto
t^{-1}t_{\{1,2,\ldots,m \}J_{a \alpha}}^{\wedge m},
\end{equation}
to be described below.

Consider the embedding of $U_q \mathfrak{sl}_{N}$-module algebras
$\mathcal{I}:\mathbb{C}[\mathcal{M}at_{m,n}]_q \hookrightarrow
\mathbb{C}[SL_N]_{q,t}$ (see the previous section) and a similar embedding
$\overline{\mathcal{I}}:\mathbb{C}[\overline{\mathcal{M}at}_{m,n}]_q
\hookrightarrow \mathbb{C}[SL_N]_{q,t^*}$
$$
\overline{\mathcal{I}}f=(\mathcal{I}f^*)^*,\qquad f \in
\mathbb{C}[\overline{\mathcal{M}at}_{m,n}]_q.
$$
(We use the obvious embeddings of the localizations
$\mathbb{C}[SL_N]_{q,t}\subset \mathbb{C}[\widetilde{G}]_{q,x}$,
$\mathbb{C}[SL_N]_{q,t^*}\subset \mathbb{C}[\widetilde{G}]_{q,x}$.) Consider
the linear map $\mathrm{Pol}(\mathcal{M}at_{m,n})_q \to
\mathbb{C}[\widetilde{G}]_{q,x}$, $f_- \cdot f_+ \mapsto
\overline{\mathcal{I}}(f_-) \cdot \mathcal{I}(f_+)$, $f_-\in
\mathbb{C}[\overline{\mathcal{M}at}_{m,n}]_q$, $f_+\in
\mathbb{C}[\mathcal{M}at_{m,n}]_q$. By our construction, this map is a
morphism of $U_q \mathfrak{su}_{n,m}$-modules and satisfies (\ref{cI}).
Prove that it is a homomorphism of $*$-algebras. It suffices to show that
$$
(\mathcal{I}z_b^\beta)^*(\mathcal{I}z_a^\alpha)=m \sigma
R_{\overline{\mathcal{I}}\mathbb{C}[\overline{\mathcal{M}at}_{m,n}]_q
\:\mathcal{I}\mathbb{C}[\mathcal{M}at_{m,n}]_q}
((\overline{\mathcal{I}}z_b^\beta)^*\otimes(\mathcal{I}z_a^\alpha))
$$
for all $a,b=1,2,\ldots,n$; $\alpha,\beta=1,2,\ldots,m$. For that, it
suffices to establish
$$
\overline{\mathcal{I}}((z_b^\beta)^*)t^{*j}\cdot t^k
\mathcal{I}(z_a^\alpha)=q^{{\tt const}\cdot j \cdot k}m \sigma
R_{\mathbb{C}[SL_N]_{q,t^*}\mathbb{C}[SL_N]_{q,t}}
(\overline{\mathcal{I}}((z_b^\beta)^*)t^{*j}\otimes t^k
\mathcal{I}(z_a^\alpha)),
$$
with $j,k \in \mathbb{Z}$. We may restrict ourselves to the special case
$j,k \in \mathbb{N}$ since
$$
t^{-k}(q^{{\tt const}\cdot j \cdot k}m \sigma
R_{\mathbb{C}[SL_N]_{q,t^*}\mathbb{C}[SL_N]_{q,t}}
(\overline{\mathcal{I}}((z_b^\beta)^*)t^{*j}\otimes t^k
\mathcal{I}(z_a^\alpha)))t^{*-j},
$$
viewed as a function of the parameter $q^{1/s}$, is a Laurent polynomial.
 In the above special case
$\overline{\mathcal{I}}((z_b^\beta)^*)t^{*j}$, $ t^k
\mathcal{I}(z_a^\alpha)$, $a,b=1,2,\ldots,n$, $\alpha,\beta=1,2,\ldots,m$,
are the matrix elements of finite dimensional weighted representations of
$U_q \mathfrak{sl}_N$. What remains is to apply the well known \cite{CP}
R-matrix commutation relations between those matrix elements, (\ref{rmmf}),
and the relations
$$
\left \langle t^k \mathcal{I}(z_a^\alpha),F_i \xi \right \rangle=\left
\langle \overline{\mathcal{I}}((z_b^\beta)^*) t^{*j},E_i \xi \right
\rangle=0,
$$
for all $\xi \in U_q \mathfrak{sl}_N$, $a,b=1,2,\ldots,n$,
$\alpha,\beta=1,2,\ldots,m, i=1,2,\ldots,N-1$. \hfill $\square$

\medskip

{\sc Remark.} A great deal of our techniques can be transferred onto the
general case of q-Cartan domains considered in \cite{SV2}. The only
exception here constitute the proofs of proposition \ref{emb1} and lemma
\ref{i'i''}. In the passage to general q-Cartan domains, the displacement of
lines in the matrix $\mathbf{t}$ used in these proofs is irrelevant.
Instead, one should use a q-analogue for the longest element
$\widetilde{w}_0$ of the Weyl group \cite{CP} such that
$\triangle(\widetilde{w}_0)=\widetilde{w}_0 \otimes \widetilde{w}_0 \cdot
R$, with $R$ being the universal R-matrix (cf. \cite[section 16]{SV2}).

\bigskip

\section{\boldmath $U_q \mathfrak{su}_{n,m}$-module algebra
$\mathrm{Pol}(\mathrm{Mat}_{m,n})_q$}\label{pol}

It was demonstrated in section \ref{local} that the structure of $U_q
\mathfrak{sl}_N$-module algebra in $\mathbb{C}[\mathcal{M}at_{m,n}]_q$ is
determined by the relations formulated in corollary \ref{uqslmat}. Our
immediate intention is to use these formulae to define a $U_q
\mathfrak{sl}_N$-module algebra structure in
$\mathbb{C}[\mathrm{Mat}_{m,n}]_q$. If $q\in(0\,,\,1)$ is a transcendental
number then the result of proposition \ref{isomo} is still valid, and we
transfer the $U_q \mathfrak{sl}_N$-module algebra structure to
$\mathbb{C}[\mathrm{Mat}_{m,n}]_q$ via the isomorphism
$$
j:\mathbb{C}[\mathrm{Mat}_{m,n}]_q\to \mathbb{C}[\mathcal{M}at_{m,n}]_q$$
defined in (\ref{zz}).
Note that this $U_q \mathfrak{sl}_N$-module algebra structure in
$\mathbb{C}[\mathrm{Mat}_{m,n}]_q$ is well defined for arbitrary $q$ as well
since the coefficients in decompositions of the elements
$z_{a_1}^{\alpha_1}z_{a_2}^{\alpha_2}\cdots z_{a_M}^{\alpha_M}$ and
$\xi(z_{a_1}^{\alpha_1}z_{a_2}^{\alpha_2}\cdots z_{a_M}^{\alpha_M})$,
$\xi\in U_q \mathfrak{sl}_N$, in the basis of lexicographically ordered
monomials \footnote{It was noted in section \ref{smr} that the linear
independence of these monomials can be proved by an application of the
Bergman diamond lemma.}
$$
(z_1^1)^{j^1_1}(z_2^1)^{j^1_2}
\cdots(z_{n-1}^m)^{j_{n-1}^m}(z_n^m)^{j_n^m},\qquad j_a^\alpha \in
\mathbb{Z}_+
$$
are Laurent polynomials in the parameter $q^{1/2}$.

Thus, $\mathbb{C}[\mathrm{Mat}_{m,n}]_q$ becomes a $U_q
\mathfrak{sl}_N$-module algebra, and we still have the homomorphism of $U_q
\mathfrak{sl}_N$-module algebras
\begin{equation}\label{zzz}
j:\mathbb{C}[\mathrm{Mat}_{m,n}]_q\to
\mathbb{C}[\mathcal{M}at_{m,n}]_q,\qquad j:z_a^\alpha \mapsto
z_a^\alpha,\quad a=1,\ldots,n,\;\alpha=1,\ldots,m.
\end{equation}
The result of proposition \ref{isomo} turns out to hold in the case of
arbitrary $q\in(0\,,\,1)$:
\medskip
\begin{proposition}\label{isom}
The homomorphism (\ref{zzz}) is an isomorphism.
\end{proposition}
{\bf Proof.} Since $j$ respects the gradation and
$$
 \dim
\mathbb{C}[\mathrm{Mat}_{m,n}]_{q,i}=\dim
\mathbb{C}[\mathcal{M}at_{m,n}]_{q,i}=\dbinom{mn+i-1}{i}
$$
(see section \ref{smr} and the equality (\ref{dimm})), it suffices to prove
injectivity of $j$.

Suppose the kernel of $j$ is a non-trivial ideal $J$. Clearly, $J$ is a $U_q
\mathfrak{sl}_N$-submodule of $\mathbb{C}[\mathrm{Mat}_{m,n}]_q$ since $j$
is a morphism of $U_q \mathfrak{sl}_N$-modules.
$\mathbb{C}[\mathrm{Mat}_{m,n}]_q$ is a $U_q\mathfrak{b}^-$-locally finite
dimensional weight $U_q\mathfrak{sl}_{N}$-module, hence the same is true for
$J$. In particular, $J$ contains a non-zero element $f$ satisfying the
relations $$
 F_if=0,\qquad i=1,\ldots,N-1.
$$
The action of $K^{\pm1}_i$ respects the above equations, hence we may assume
that $f$ is a {\it weight} vector in the $U_q\mathfrak{sl}_{N}$-module
$\mathbb{C}[\mathrm{Mat}_{m,n}]_q$, that is a common eigenvector of all
$K^{\pm1}_i$, $i=1,\ldots,N-1$.

\medskip
\begin{lemma}
 Any weight vector in
$\mathbb{C}[\mathrm{Mat}_{m,n}]_q$, annihilated by all $F_i$, $i\neq n$,
coincides up to a constant with one of the vectors
$$
f_{k_1,\ldots,k_m}=(z_n^m)^{k_1}\left(z_{\quad \{n-1,n\}}^{\wedge
2\{m-1,m\}}\right)^{k_{2}}\left(z_{\quad \{n-2,n-1,n\}}^{\wedge
3\{m-2,m-1,m\}}\right)^{k_{3}}\cdots \left(z_{\quad
\{n-m+1,\ldots,n\}}^{\wedge m\{1,\ldots,m\}}\right)^{k_{m}}
$$
($k_1,\ldots,k_m\in\mathbb{Z}_+$).
\end{lemma}
{\bf Proof} of the lemma.
 Using the explicit formulae, given in corollary
\ref{uqslmat}, one proves that the elements $f_{k_1,\ldots,k_m}$ are weight
vectors, annihilated by all $F_i$ with $i\neq n$. Impose the notation
$L_{k_1,\ldots,k_m}$ for the finite dimensional simple $U_q \mathfrak{sl}_n
\otimes U_q \mathfrak{sl}_m$-submodule in $\mathbb{C}[\mathrm{Mat}_{m,n}]_q$
generated by the vector $f_{k_1,\ldots,k_m}$. Clearly,
\begin{equation}\label{inclu}
\mathbb{C}[\mathrm{Mat}_{m,n}]_{q,i}\supset\bigoplus L_{k_1,\ldots,k_m}
\end{equation}
where the sum is taken over all $m$-tuples $(k_1,\ldots,k_m)$ satisfying
$k_1+2k_2+\ldots+mk_m=i$. In the classical case one has the equality in
(\ref{inclu}). Thus one has the equality in the quantum case as well since
the dimensions of the spaces $\mathbb{C}[\mathrm{Mat}_{m,n}]_{q,i}$,
$L_{k_1,\ldots,k_m}$ are just the same as in the classical case. We conclude
that any homogeneous (in particular, weight) vector, annihilated by all
$F_i$ with $i\neq n$, is a linear combination of some $f_{k_1,\ldots,k_m}$'s
as the lowest weight vector in a simple $U_q \mathfrak{sl}_n \otimes U_q
\mathfrak{sl}_m$-module is unique up to a constant. To finish the proof of
the lemma it remains to observe that weights of the vectors
$f_{k_1,\ldots,k_m}$ are distinct for different $m$-tuples
$(k_1,\ldots,k_m)$. \hfill$\square$

\medskip

Return to the proof of proposition \ref{isom}. Due to the above lemma we may
assume that $J$ contains $f_{k_1,\ldots,k_m}$ for some
$k_1,\ldots,k_m\in\mathbb{Z}_+$, $k_1+\ldots+k_m\neq0$, which, moreover,
satisfies $F_nf_{k_1,\ldots,k_m}=0$.
 Let
$$ f_{k_1,\ldots,k_m}=
\sum_{j=0}^s\psi_{j\,;\,k_1,\ldots,k_m}\cdot(z_n^m)^{j},\quad
s\in\mathbb{Z}_+$$ where $\psi_{j\,;\,k_1,\ldots,k_m}$ are elements of the
unital subalgebra in $\mathbb{C}[\mathrm{Mat}_{m,n}]_q$ generated by
$\{z_a^\alpha\}_{(\alpha,a)\neq(m,n)}$. Then, by corollary \ref{uqslmat}
$$
F_nf_{k_1,\ldots,k_m}=q^{1/2}\cdot
\sum_{j=0}^s\frac{1-q^{-2j}}{1-q^{-2}}\cdot\psi_{j\,;\,k_1,\ldots,k_m}\cdot(z_n^m)^{j-1}.
$$
So $F_nf_{k_1,\ldots,k_m}=0$ implies
$\psi_{k_1+k_2+\ldots+k_m\,;\,k_1,\ldots,k_m}=0$ because of
$k_1+k_2+\ldots+k_m\neq0$. On the other hand,
$\psi_{k_1+k_2+\ldots+k_m\,;\,k_1,\ldots,k_m}\neq0$ since the image of
$f_{k_1,\ldots,k_m}$ under the homomorphism of algebras
$$\varphi:\mathbb{C}[\mathrm{Mat}_{m,n}]_q\to\mathbb{C}[z],\quad
\varphi:z_a^\alpha\mapsto\begin{cases} z,&a=n \;\&\;\alpha=m,
\\ 0,&a-\alpha \ne n-m,
\\ 1,&\mathrm{otherwise},
\end{cases}$$
is $z^{k_1+k_2+\ldots+k_m}\neq0$.
 \hfill$\square$

\medskip

Let us endow $\mathbb{C}[\overline{\mathrm{Mat}}_{m,n}]_q$ with a $U_q
\mathfrak{sl}_N$-module algebra structure via
$$
(\xi f)^*=(S(\xi))^*f^*,\qquad \xi \in U_q \mathfrak{su}_{n,m},\quad f \in
\mathbb{C}[\overline{\mathrm{Mat}}_{m,n}]_q.
$$
Proposition \ref{isom} implies that there is a canonical isomorphism of $U_q
\mathfrak{sl}_N$-module algebras
\begin{equation}\label{zzz*}
\mathbb{C}[\overline{\mathrm{Mat}}_{m,n}]_q\to
\mathbb{C}[\overline{\mathcal{M}at}_{m,n}]_q,\quad (z_a^\alpha)^* \mapsto
(z_a^\alpha)^*,\quad a=1,\ldots,n,\;\alpha=1,\ldots,m.
\end{equation}

\medskip

{\sc Remark.} Note that for any $q\in(0,1)$ the $*$-algebra
$\mathrm{Pol}(\mathcal{M}at_{m,n})_q$ is generated by $z_a^\alpha$,
$a=1,\ldots,n$, $\alpha=1,\ldots,m$ due to the isomorphisms (\ref{zzz}),
(\ref{zzz*}). The same fact for transcendental $q$ is an easy consequence of
results of section \ref{da}.

\medskip

\begin{proposition}\label{ptp}\hfill \\
i) The homomorphism (\ref{isomor}) is an isomorphism of $*$-algebras;

\noindent ii) there exists a unique structure of $U_q
\mathfrak{su}_{n,m}$-module algebra in $\mathrm{Pol}(\mathrm{Mat}_{m,n})_q$
such that (\ref{isomor}) is an isomorphism of $U_q
\mathfrak{su}_{n,m}$-module algebras.
\end{proposition}

{\bf Proof.} Statement ii) immediately follows from i).

To prove i), consider the filtration
$$F^k\mathrm{Pol}(\mathrm{Mat}_{m,n})_q=\bigoplus_{i+j\leq k}
\mathrm{Pol}(\mathrm{Mat}_{m,n})_{q,i,-j}$$ (see (\ref{bigrad})) and a
similar filtration $F^k\mathrm{Pol}(\mathcal{M}at_{m,n})_q$,
$k\in\mathbb{Z}_+$. The homomorphism (\ref{isomor}) respects these
filtrations. Clearly,

$$\dim F^k\mathrm{Pol}(\mathrm{Mat}_{m,n})_q\leq\sum_{j=0}^k
\dbinom{2mn+j-1}{j}.$$ On the other hand, due to the isomorphisms
(\ref{zzz}), (\ref{zzz*})
$$\dim F^k\mathrm{Pol}(\mathcal{M}at_{m,n})_q=\sum_{j=0}^k
\dbinom{2mn+j-1}{j}.$$ It remains to observe that the homomorphism
(\ref{isomor}) is surjective as $z_a^\alpha$, $a=1,\ldots,n$,
$\alpha=1,\ldots,m$, generate the $*$-algebra
$\mathrm{Pol}(\mathcal{M}at_{m,n})_q$.\hfill$\square$

\medskip
\begin{corollary}\label{corol}
$
\mathbb{C}[\mathrm{Mat}_{m,n}]_q \otimes
\mathbb{C}[\overline{\mathrm{Mat}}_{m,n}]_q \simeq
\mathrm{Pol}(\mathrm{Mat}_{m,n})_q.
$
\end{corollary}

\bigskip

\section{A proof of proposition \ref{yexp}}\label{pry}

It follows from the definitions of section \ref{pol} and the remark before
lemma \ref{Noumi} that the homomorphism $i:\mathbb{C}[\mathrm{Mat}_{m,n}]_q
\to \mathbb{C}[SL_N]_{q,t}$ is a morphism of $U_q \mathfrak{sl}_N$-module
algebras.

\smallskip

\begin{proposition}\label{cqk}
Let $k \in{\mathbb N}$. There exists a constant $c(q,k)\in \mathbb{C}$ such
that for all $1 \le \alpha_1<\alpha_2<\dots<\alpha_k \le m$, $1 \le
a_1<a_2<\dots<a_k \le n$, in the algebra $\mathbb{C}[SL_N]_{q,t}$, one has
\begin{equation}\label{Iz}
i \left(z_{\quad \{a_1,a_2,\dots,a_k \}}^{\wedge k
\{m+1-\alpha_k,m+1-\alpha_{k-1},\dots,m+1-\alpha_1
\}}\right)=c(q,k)t^{-1}t_{\{1,2,\dots,m \}J}^{\wedge m},
\end{equation}
with $J=\{n+1,n+2,\dots,N \}\setminus
\{n+\alpha_1,n+\alpha_2,\dots,n+\alpha_k \}\cup \{a_1,a_2,\dots,a_k \}$.
\end{proposition}

\smallskip

{\bf Proof.} An argument similar to that used while proving lemma
\ref{Noumi} reduces the general case to the special case $m=n=k$:
\begin{equation}\label{Iz1}
i \left(z_{\quad \{1,2,\dots,k \}}^{\wedge k \{1,2,\dots,k
\}}\right)=c(q,k)t^{-1}t_{\{1,2,\dots,k \}\{1,2,\dots,k \}}^{\wedge k}.
\end{equation}
Let $\widetilde{z}=z_{\quad \{1,2,\dots,k \}}^{\wedge k \{1,2,\dots,k \}}$,
$\widetilde{t}=t_{\{1,2,\dots,k \}\{1,2,\dots,k \}}^{\wedge k}$, and $U_q
\mathfrak{sl}_k \otimes U_q \mathfrak{sl}_k \subset U_q \mathfrak{sl}_{2k}$
the Hopf subalgebra generated by $E_j$, $F_j$, $K_j^{\pm 1}$, $j \ne k$.

Consider the subalgebra $\mathcal{F}$ of $\mathbb{C}[SL_{2k}]_q$ generated
by the quantum minors $t_{\{1,2,\dots,k \}J}^{\wedge k}$,
$\mathrm{card}\,J=k$. We need the following

\begin{lemma}\label{bas}
$\left \{\left.t^{j'}\widetilde{t}^{j''}\right|\,j',j''\in
\mathbb{Z}_+\right \}$ form a basis in the vector space of $U_q
\mathfrak{sl}_k \otimes U_q \mathfrak{sl}_k$-invariants in $\mathcal{F}$.
\end{lemma}

{\bf Proof.} It is obvious that $t^{j'}\widetilde{t}^{j''}$ are $U_q
\mathfrak{sl}_k \otimes U_q \mathfrak{sl}_k$-invariant. They are linear
independent since for all $\xi \in U_q \mathfrak{sl}_{2k}$
\begin{equation}\label{weights}
\left \langle t^{j'}\widetilde{t}^{j''},K_k \xi \right
\rangle=q^{j'+j''}\left \langle t^{j'}\widetilde{t}^{j''},\xi \right
\rangle,\qquad \left \langle t^{j'}\widetilde{t}^{j''},\xi K_k \right
\rangle=q^{j''-j'}\left \langle t^{j'}\widetilde{t}^{j''},\xi \right
\rangle.
\end{equation}
Prove that the above vectors generate the space of $U_q \mathfrak{sl}_k
\otimes U_q \mathfrak{sl}_k$-invariants. $\left \{\left.\widetilde{t}^l
\right|\,l \in \mathbb{Z}_+\right \}$ form a weight basis of the vector
space $\{f \in \mathcal{F}|\;E_jf=0,\:j=1,2,\ldots,N-1\}$. Hence the $U_q
\mathfrak{sl}_{2k}$-module $\mathcal{F}$ is isomorphic to a sum of the
simple $U_q \mathfrak{sl}_{2k}$-modules $U_q
\mathfrak{sl}_{2k}\widetilde{t}^l$, $l \in \mathbb{Z}_+$. The dimensions of
their weight subspaces remain intact under the passage from the classical to
the quantum case (see \cite{Jant}). Therefore the dimensions of their $U_q
\mathfrak{sl}_k \otimes U_q \mathfrak{sl}_k$-isotypic components also remain
intact. Thus we conclude that $\left
\{\left.t^{j'}\widetilde{t}^{j''}\right|\,j'+j''=l \right \}$ generate the
vector space of $U_q \mathfrak{sl}_k \otimes U_q \mathfrak{sl}_k$-invariants
in $U_q \mathfrak{sl}_{2k}\widetilde{t}^l$ and $\left
\{\left.t^{j'}\widetilde{t}^{j''}\right|\,j',j''\in \mathbb{Z}_+\right \}$
generate the vector space of $U_q \mathfrak{sl}_k \otimes U_q
\mathfrak{sl}_k$-invariants in $\mathcal{F}$. \hfill $\square$

\medskip

Turn back to the proof of proposition \ref{cqk} in the special case $m=n=k$.
Observe that $i(\widetilde{z})$ belongs to the localization of $\mathcal{F}$
with respect to the multiplicative set $t^\mathbb{N}$ and is $U_q
\mathfrak{sl}_N$-invariant. Hence by lemma \ref{bas} $i(\widetilde{z})$
belongs to the linear span of $\left
\{\left.t^{j'}\widetilde{t}^{j''}\right|\,j'\in \mathbb{Z},j''\in
\mathbb{Z}_+\right \}$. What remains is to compare (\ref{weights}) and
$$
\left \langle i(\widetilde{z}),K_k \xi \right \rangle=0,\qquad \left \langle
i(\widetilde{z}),\xi K_k \right \rangle=q^2 \left \langle
i(\widetilde{z}),\xi \right \rangle. \eqno \square
$$

\medskip

{\bf Proof} of lemma \ref{Izm}. To compute the constant $c(q,k)$ in
(\ref{Iz}), apply $F_n$ to both sides of this relation. The relations
$F_nz_a^\alpha=q^{1/2}\delta_{an}\delta_{\alpha m}$ and $H_nz_b^\beta=0$ for
$b \ne n$, $\beta \ne m$ (see corollary \ref{uqslmat}) imply
$$
F_nz_{\quad \{n-k+1,\dots,n \}}^{\wedge k \{m-k+1,\dots,m
\}}=q^{1/2}z_{\quad \{n-k+1,\dots,n-1 \}}^{\wedge k-1 \{m-k+1,\dots,m-1
\}}.
$$
It follows from $\Delta(F_n)=F_n \otimes K_n^{-1}+1 \otimes F_n$,
$F_n(t^{-1})=0$ that
\begin{multline*}
F_n \left(t^{-1}t_{\{1,2,\dots,m \}\{n-k+1,\dots,n,n+k+1,\dots,N \}}^{\wedge
m}\right)=t^{-1}F_nt_{\{1,2,\dots,m \}\{n-k+1,\dots,n,n+k+1,\dots,N
\}}^{\wedge m}=
\\ =q^{1/2}t^{-1}t_{\{1,2,\dots,m \}\{n-k+1,\dots,n-1,n+1,n+k+1,\dots,N
\}}^{\wedge m}.
\end{multline*}
Hence one has in (\ref{Iz}) $c(q,k)=c(q,k-1)=\dots=c(q,1)=1$. \hfill
$\square$

\medskip

\begin{lemma}
Let ${\rm card}(J)=m$, $J^c=\{1,2,\dots,N \}\setminus J$,
$l(J,J^c)=\mathrm{card}\{(j',j'')\in J \times J^c|\;j'>j'' \}$. Then
\begin{equation}\label{min*}
\left(t_{\{1,2,\dots,m \}J}^{\wedge
m}\right)^*=(-1)^{\mathrm{card}(\{1,2,\dots,n \}\cap
J)}(-q)^{l(J,J^c)}t_{\{m+1,m+2,\dots,N \}J^c}^{\wedge n}.
\end{equation}
\end{lemma}

\smallskip

{\bf Proof.} This lemma is deducible from (\ref{star}) and a general formula
of Ya. Soibelman \cite{So1}, \cite[p. 432]{CP}:
$$
\left(t_{\{1,2,\dots,m \}J}^{\wedge
m}\right)^\star=(-q)^{l(J,J^c)}t_{\{m+1,m+2,\dots,N \}J^c}^{\wedge n}.\eqno
\square
$$

\medskip

It now follows from (\ref{qdet}), (\ref{min*}) that
\begin{equation}\label{tmtm*}
\sum_{\genfrac{}{}{0pt}{1}{J \subset \{1,\dots,N
\}}{\mathrm{card}(J)=m}}(-1)^{\mathrm{card}(\{1,2,\dots,n \}\cap
J)}t_{\{1,2,\dots,m \}J}^{\wedge m}\left(t_{\{1,2,\dots,m \}J}^{\wedge
m}\right)^*=1.
\end{equation}

\smallskip

{\bf Proof} of proposition \ref{yexp}. Multiply (\ref{tmtm*}) by $t^{-1}$
from the left and by $t^{* -1}$ from the right and apply lemma \ref{Izm},
together with (\ref{tmtm*}). The result is just the statement of proposition
\ref{yexp}. \hfill $\square$

\bigskip

\section{Faithfulness of the representation \boldmath $T$}\label{ii}

We start with proving the faithfulness of the representation $T$.

\begin{lemma}\label{oo}
For all $k,l \in \mathbb{Z}_+$, the map
$$
\mathbb{C}[\mathrm{Mat}_{m,n}]_{q,k}\otimes
\mathbb{C}[\overline{\mathrm{Mat}}_{m,n}]_{q,-l}\to
\mathrm{Hom}(\mathcal{H}_l,\mathcal{H}_k);\qquad f_+ \otimes f_-\mapsto
T(f_+f_-)|_{\mathcal{H}_l}
$$
is one-to-one.
\end{lemma}

\smallskip

{\bf Proof.} It follows from corollary \ref{sh0} that
$$
i_1:\mathbb{C}[\mathrm{Mat}_{m,n}]_{q,k}\to \mathcal{H}_k,\qquad
i_1:f_+\mapsto f_+v_0,\qquad f_+\in \mathbb{C}[\mathrm{Mat}_{m,n}]_{q,k},
$$
is an isomorphism of vector spaces. Since $\mathcal{H}\simeq \mathscr{L}$
(see proposition \ref{iso}) and $\mathcal{H}$ is a pre-Hilbert space, one
has an isomorphism
$$
i_2:\mathbb{C}[\overline{\mathrm{Mat}}_{m,n}]_{q,-l}\to
\mathcal{H}_l^*,\qquad \left \langle i_2(f_-),v \right
\rangle=(v,f_-^*v_0),\qquad f_-\in
\mathbb{C}[\overline{\mathrm{Mat}}_{m,n}]_{q,-l}.
$$
What remains is to use the canonical isomorphism $\mathcal{H}_k \otimes
\mathcal{H}_l^*\simeq \mathrm{Hom}(\mathcal{H}_l,\mathcal{H}_k)$ and the
relation
$$
T(f_+f_-)v=T(f_+)(T(f_-)v,v_0)v_0=(v,f_-^*v_0)f_+v_0
$$
for $v \in \mathcal{H}_l$, $f_+\in \mathbb{C}[\mathrm{Mat}_{m,n}]_{q,k}$,
$f_-\in \mathbb{C}[\overline{\mathrm{Mat}}_{m,n}]_{q,-l}$.\hfill $\square$

\medskip

\begin{proposition}\label{em}
The representation $T$ is faithful.
\end{proposition}

\smallskip

{\bf Proof.} Recall that the vector space
$\mathrm{Pol}(\mathrm{Mat}_{m,n})_q$ is equipped with a bigradation
(\ref{bigrad}). We need a standard partial order on $\mathbb{Z}_+^2$:
$$
(k_1,l_1)\le(k_2,l_2)\qquad \Leftrightarrow \qquad k_1 \le k_2 \quad \&
\quad l_1 \le l_2.
$$

Assume that our statement is wrong and $T(f)=0$ for some $f \in
\mathrm{Pol}(\mathrm{Mat}_{m,n})_q$, $f \ne 0$. Consider a homogeneous
component $f_{kl}\ne 0$ of $f$ with minimal bidegree $(k,l)$. (Such
homogeneous component certainly exists, but it may be non-unique for a given
$f \in \mathrm{Pol}(\mathrm{Mat}_{m,n})_q$). Let $P_k:\mathcal{H}\to
\mathcal{H}_k$ be the projection onto $\mathcal{H}_k$ with kernel $\bigoplus
\limits_{j \ne k}\mathcal{H}_j$. Since $f_{kl}$ is of a minimal bidegree,
one has $P_kT(f_{kl})|_{\mathcal{H}_l}=P_kT(f)|_{\mathcal{H}_l}=~0$,
$f_{kl}\ne 0$, which contradicts the statement of lemma \ref{oo}. \hfill
$\square$

\medskip

{\sc Remark.} Recall that in section \ref{local} we constructed a
homomorphism $i:\mathrm{Pol}(\mathrm{Mat}_{m,n})_q\to
\mathbb{C}[\widetilde{G}]_{q,x}$ (see proposition \ref{emb3}). Now we may
establish its injectivity. It follows from proposition \ref{em} and the
commutative diagram:
$$
\begin{CD}
\mathrm{Pol}(\mathrm{Mat}_{m,n})_q @>i>> \mathbb{C}[\widetilde{G}]_{q,x}
\\ @VTVV @VV{\widetilde{\mathcal{T}}}V
\\ \mathrm{End}(\mathcal{H})@>\thicksim>> \mathrm{End}(\mathscr{L})
\end{CD}
$$

\pagebreak

\end{document}